\documentclass{amsart}
\usepackage{amssymb,mathrsfs,stmaryrd,imakeidx,hyperref,enumitem}
\input xy \xyoption {all}
\usepackage{tikz}
\usetikzlibrary{cd}
\usepackage{dsfont}
\usepackage{tikz-cd}
\usepackage{epigraph}
\makeindex

\numberwithin{equation}{section}
\numberwithin{figure}{section}
\newtheorem{thm}{Theorem}[section]
\newtheorem{lem}[thm]{Lemma}
\newtheorem{prop}[thm]{Proposition}
\newtheorem{cor}[thm]{Corollary}

\theoremstyle{definition}
\newtheorem{defn}[thm]{Definition}
\newtheorem{defn-prop}[thm]{Definition-Proposition}

\theoremstyle{remark}
\newtheorem{rmk}[thm]{Remark}
\newtheorem*{Sketch*}{Sketch of proof}
\newtheorem{ex}[thm]{Example}
\usepackage{url}

\usepackage{abstract}

\usepackage{amsmath,latexsym}
\usepackage{mhchem}
\usepackage{xcolor}
\usepackage[utf8x]{inputenc}
\usepackage{extarrows}
\usepackage{chemarrow}
\graphicspath{{figure/}}                                 
\usepackage[a4paper]{geometry}                           
\newcommand{\nc}{\newcommand}

\nc\on{\operatorname}
\nc\op{\operatorname}

\renewcommand{\c}[1]{\mathcal{#1}}

\renewcommand{\b}[1]{\mathbb{#1}}

\newcommand{\f}[1]{\mathfrak{#1}}

\renewcommand{\r}[1]{\mathrm{#1}}

\nc{\h}{{\hslash}}	

\nc{\ul}[1]{\underline{#1}}
\nc{\ol}[1]{\overline{#1}}
\nc{\wh}[1]{\widehat{#1}}
\nc{\wt}[1]{\widetilde{#1}}

\nc{\bs}{\backslash}

\nc\sm{\setminus}

\nc{\Aut}{\op{Aut}}
\nc{\Tor}{\op{Tor}}
\nc{\Ext}{\op{Ext}}
\nc{\DMod}{\op{DMod}}
\nc{\Span}{\op{Span}}
\nc{\gr}{\op{gr}}
\nc{\Vect}{\op{Vect}}
\nc\ch{\op{ch}}
\nc\Rep{\op{Rep}}
\nc\diag{\op{diag}}
\nc\Tr{\op{Tr}}
\nc\tr{\op{tr}}
\nc\codim{\op{codim}}
\nc\Sym{\op{Sym}}
\nc\length{\op{length}}
\nc\Stab{\op{Stab}}
\nc\RHom{\op{RHom}}
\nc\Alt{\op{Alt}} 

\renewcommand{\sl}{\mathfrak{sl}}

\newcommand{\SL}{\r{SL}}

\nc{\p}{\partial}
\nc{\Lie}{\op{Lie}}
\nc{\Ug}{\mathcal{U}\mathfrak{g}}      
\nc{\UX}{\cU_{\cO_X}} 
\nc{\ad}{\op{ad}}
\nc{\Ad}{\op{Ad}}
\nc{\height}{\op{ht}}

\nc{\inj}{\hookrightarrow}
\nc{\surj}{\twoheadrightarrow}
\nc\Ind{\on{Ind}}
\nc\Res{\on{Res}}
\nc\Loc{\on{Loc}}
\nc\Map{\op{Map}}
\nc\Qmap{\op{QMap}}
\nc\Hom{\op{Hom}}
\nc{\Mod}{\op{Mod}}
\nc\End{\op{End}}
\nc\im{\op{im}}
\nc\coker{\op{coker}}
\nc\Dis{\op{Dis}} 
\nc\Isom{\op{Isom}}
\nc\pr{\op{pr}}
\nc\ladic{\ell\op{-adic}}
\nc\padic{p\op{-adic}}
\nc\Gal{\op{Gal}}

\nc\Pic{\op{Pic}}
\nc\Spec{\op{Spec}}
\nc{\Sh}{\mathsf{Sh}}
\nc\QCoh{\mathsf{QCoh}}
\nc\Bun{\op{Bun}}
\nc{\IC}{\on{IC}}
\nc{\Perv}{\mathsf{Perv}}
\nc{\itquot}{/\!\!/}
\nc{\Gr}{\op{Gr}}
\nc{\Coh}{\mathsf{Coh}}
\nc{\Div}{\op{Div}}
\nc{\supp}{\op{supp}}
\nc{\Hilb}{\op{Hilb}}
\nc{\GHilb}{\Gamma\text{-}\op{Hilb}}
\nc{\Ran}{\op{Ran}}
\nc{\Picf}{\Pic^{\op{fact}}}
\nc{\Picsf}{\Pic^{\op{sf}}}
\nc{\Perf}{\op{Perf}}
\nc{\Nm}{\op{Nm}}
\nc{\Set}{\mathsf{Set}}
\nc{\FinSet}{\mathsf{FinSet}}
\nc{\Ab}{\mathsf{Ab}}
\nc{\Cat}{\mathsf{Cat}}
\nc{\Stack}{\mathsf{Stack}}
\nc{\Top}{\mathsf{Top}}
\nc{\Ob}{\op{Ob}}
\nc{\Mor}{\op{Mor}}
\nc{\id}{\op{id}}
\nc{\Sch}{\mathsf{Sch}}
\nc{\Indsch}{\mathsf{IndSch}}
\nc{\Fun}{\op{Funct}}
\nc{\corre}[5]{\left(#1\xleftarrow{#2}#3\xrightarrow{#4}#5\right)}

\nc{\Sol}{\op{Sol}}
\nc{\LS}{\op{LocSys}}

\nc{\GK}{G_\mathcal{K}}
\nc{\GO}{G_\mathcal{O}}
\nc{\Wa}{W_{\op{aff}}}
\nc{\Fl}{\c{F}\ell}
\nc{\Int}{\op{Int}}
\nc{\reg}{\op{reg}}

\begin{document}
	
	\title{\bf Affine Grassmannian from divisors}
\author[Chunlin Shao]{Chunlin Shao}
\address{University of Massachusetts, Amherst, MA, USA.}
\email{cshao@umass.edu}
	\maketitle
\begin{abstract}
We generalize Mirkovi\'c's construction of Zastava spaces. These spaces form an inductive system, giving a generalization of affine and Beilinson--Drinfeld Grassmannians as indschemes, defined for any integral symmetric matrix $\kappa$ with even diagonal entries. When $\kappa$ is the  Cartan matrix from a simply-laced Kac--Moody algebra, we also construct an action of the Tits group on the generalized affine Grassmannian.
	\end{abstract}
	\vspace{1ex}	
	\vspace{1ex}

\section{Introduction}
\subsection{What is this paper about?}
The present thesis develops a new construction of generalized affine
Grassmannians from Zastava spaces.
The affine Grassmannian $\Gr_G$ of a reductive algebraic group $G$
plays an important role in representation theory, number theory,
algebraic geometry, particularly in the
geometric Langlands program.
See \cite{zhu2016introduction} for an introduction.

One way to study affine Grassmannians is through their finite-dimensional
pieces. The pieces relevant to this paper arise as fibers of Zastava spaces, introduced by Finkelberg and Mirkovi\'c \cite{finkelberg1999semi} following ideas of Drinfeld.

Mirkovi\'c has extended the notion of affine Grassmannian to spaces associated with arbitrary integral
symmetric matrices $\kappa$ with even diagonal entries, including
matrices that do not arise from (reductive nor Kac--Moody) groups
\cite{mirkovic2021loop,mirkovic2023lattice}.
His approach realizes Zastava spaces by a new algebro-geometric construction called \textit{local projective spaces}, see \cite{mirkovic2014loop} for historical remarks.
The local projective space applies to any such matrices $\kappa$. Its base is the colored (i.e. a product of) Hilbert schemes of points on a smooth curve $C$, parametrizing effective colored divisors on $C$.

In this paper, we extend the constructions of Zastava spaces via
local projective spaces to pairs of colored divisors $(D,E)$,
where divisors are not required to be effective, but their
difference $E-D$ is. In such a way we define the affine Grassmannian associated to any $\kappa$ as a union of these generalized Zastavas for any pair $(D,E)$. This construction simplifies the earlier one, as functoriality is stronger and there are much fewer consistency checks.

We expect that this presentation will provide new viewpoints
for the study of affine Grassmannians, even in the classical case
associated with a reductive group.

\subsection{The classical case of a reductive group}
\label{subsec:introduction-grassmannians-zastava}

Let $G$ be a connected simple simply connected algebraic group over
$\mathbb C$. Set $\mathcal O=\mathbb C[\![z]\!]$ and
$\mathcal K=\mathbb C(\!(z)\!)$. The affine Grassmannian (also called the loop Grassmannian)
$\Gr_G:=\GK/\GO$ of $G$ is an indscheme.
Let $C$ be a smooth connected complex curve, and let
$\mathcal H_C=\coprod_{n\geq 0}C^n/S_n$ be the space of
effective divisors on $C$, i.e. the Hilbert scheme of points of $C$. Here \(S_n\) acts on \(C^n\) by permuting the factors.
Allowing the points to vary defines the
Beilinson--Drinfeld Grassmannian $\Gr^{\mathrm{BD}}_G$ over
$\mathcal H_C$, whose fiber over any point $p\in C\subset \c H_C$ is
isomorphic to $\Gr_G$. $\Gr^{\mathrm{BD}}_G$ carries a factorization
structure. See \cite{zhu2016introduction} for an introduction to
affine and Beilinson--Drinfeld Grassmannians.

Fix opposite Borel subgroups $B=TN$ and $B^-=TN^-$, and let
$I$ index the simple coroots. We identify $\mathbb Z[I]$ with
$X_*(T)$ and write $\lambda\leq\mu$ if
$\mu-\lambda\in\mathbb N[I]$.
The semi-infinite orbits in $\Gr_G$ are
\[
S_\lambda=N(\mathcal K)z^\lambda\GO/\GO,
\qquad
T_\lambda=N^-(\mathcal K)z^\lambda\GO/\GO,
\qquad \lambda\in X_*(T).
\]
Their closures satisfy
\[
\overline S_\lambda=\bigcup_{\nu\leq\lambda}S_\nu,
\qquad
\overline T_\lambda=\bigcup_{\nu\geq\lambda}T_\nu.
\]
The intersection $S_\lambda\cap\overline T_\mu$ is nonempty
precisely when $\mu\leq\lambda$
\cite{mirkovic2007geometric}.

The finite-dimensional subschemes relevant to this paper are the Zastava spaces which will be explained as follows. For a coweight $\alpha^\vee=\sum_{i\in I}a_i\alpha_i^\vee\in\mathbb N[I]$, where $\alpha_i^\vee$ are the simple coroots,
the Zastava space $\mathcal Z_G^{\alpha^\vee}$ is equipped with a
factorization morphism
\[
\pi:\mathcal Z_G^{\alpha^\vee}\longrightarrow
\mathcal H_{C\times I}^{\alpha^\vee}
:=\prod_{i\in I}C^{a_i}/S_{a_i}.
\]
When $C=\mathbb A^1$, Zastava space $\mathcal Z_G^{\alpha^\vee}$ is the space of based quasimaps of degree
$\alpha^\vee$ from $\mathbb P^1$ to the flag variety of $G$, as in
\cite{finkelberg1999semi}. 

For \(p\in C\), let \(\alpha^\vee[p]\) denote the colored divisor whose component in
color $i$ is $a_i[p]$. The fiber of $\pi$ over ${\alpha^\vee}[p]$ is identified with
\[
\mathcal Z_{G,{\alpha^\vee}[p]}
\simeq S_{0,p}\cap\overline T_{-\alpha^\vee,p}
\subset\Gr_{G,p}:=G(\c K_p)/G(\c O_p)\simeq \Gr_G,
\]
where $S_{0,p}, T_{-\alpha^\vee,p}$ denotes the corresponding semiinfinite orbits in $\Gr_{G,p}$, see \cite[\S6.4.2]{finkelberg1999semi}.
In particular, for a simple coroot $\alpha_i^\vee$, we have
$S_{0,p}\cap\overline T_{-\alpha_i^\vee,p}\simeq\mathbb A^1.$

Note that for any coweights $\lambda\le \mu$, we have $\overline S_{-\lambda}\cap\overline T_{-\mu}=z^{-\lambda}\left(\ol S_{0}\cap\overline T_{-(\mu-\lambda)}\right)$. The affine Grassmannian can be written as a direct limit in terms of Zastava fiber closures:
\[
\Gr_G\simeq
\varinjlim_{\lambda\leq\mu}
\bigl(\overline S_{-\lambda}\cap\overline T_{-\mu}\bigr).
\]
This description motivates our construction of a family of
interval Zastava spaces $Z$, indexed by pairs of colored divisors
$D\le E$, for arbitrary integral symmetric matrices $\kappa$
with even diagonal entries. 
When $\kappa$ is the Cartan matrix of a group $G$ of type ADE,
our construction recovers the intersections above: for $p\in C$, the Zastava fiber $Z_{(\lambda[p],\mu[p])}$ is canonically isomorphic to 
$\overline S_{-\lambda}\cap\overline T_{-\mu}$.

\subsection{Overview of constructions}
\subsubsection{Mirkovi\'c's Zastava construction}
\label{Mirkovic zastava construction}

For $\Gr_G$, the ample generator $\mathcal O(1)$ of $\Pic(\Gr_G)$ gives the  Kodaira embedding
$\Gr_G\hookrightarrow
\mathbb P\!\left(\Gamma(\Gr_G,\mathcal O(1))^\vee\right).$
The idea in Mirkovi\'c's work \cite{mirkovic2021loop,mirkovic2023lattice} is to reconstruct the 
Zastava spaces from their \(T\)-fixed point subschemes and the restriction of \(\mathcal O(1)\), using \((\operatorname{Gr}_G)^T=\operatorname{Gr}_T\).

Assume that $G$ is of type ADE. Write
an effective colored divisor as
$D=\sum_{p\in C} \alpha_{p}[p]$, with
$\alpha_p\in\mathbb N[I]$, and consider the closure of ordinary zastava fiber which is explained in section \ref{subsec:introduction-grassmannians-zastava}, still denoted by $Z$:
\[
 Z_{G,D}:=
\prod_{p}
\left(\overline S_{0,p}\cap\overline T_{-\alpha_p,p}\right)
\subset
\prod_{p}\Gr_{G,p}.
\]
The Beilinson--Drinfeld version $\ol S_{0,E}^{BD}$ of the semiinfinite orbit $\ol{S_0}$ for $E\in \c H_C$ is then a colimit of zastava fibers $Z_{G,D}$ over $D\subset \wh E\times I$, $\Gr_G^{BD}$ is then a union of isomorphic copies  $\ol S_\lambda^{BD}\simeq \ol S_0^{BD}$.

The closure relations for
semiinfinite orbits show that the $T$-fixed points of $ Z_{G,D}$ are
indexed by the effective colored divisors $0\leq D'\leq D$ via the Abel--Jacobi map, see \cite[\S2.2.3]{mirkovic2021loop}. In particular, the $T$-fixed point subscheme $\ol{S}_{0,p}^T$ is the Hilbert scheme $\c H_{\wh p\times I}\subset \c H_{C\times I}$, where $\wh p$ is the formal neighborhood of $p$.
To obtain their scheme structure, one considers the Hilbert scheme
$\mathcal H_D$ of all closed subschemes (of all possible lengths) of the finite colored scheme $D$, nonreduced in general. In other words, for $D=(D_i)_{i\in I}$, $\mathcal H_D=\prod_{i\in I}\mathcal H_{D_i}\subset \mathcal H_{C\times I}$.

The restriction theorem \cite[\S2.4.6]{mirkovic2021loop} gives an isomorphism
\[
\Gamma(Z_{G,D},\mathcal O(1))
\xrightarrow{\sim}
\Gamma\bigl(\mathcal H_D,\mathcal O(1)|_{\mathcal H_D}\bigr).
\]
Together with the Kodaira embedding of $\mathcal H_D$, this yields
\[
\mathcal H_D
\hookrightarrow Z_{G,D}
\hookrightarrow
\mathbb P\!\left(
\Gamma\bigl(\mathcal H_D,\mathcal O(1)|_{\mathcal H_D}\bigr)^\vee
\right).
\]
It remains to describe the image of $Z_{G,D}$ in this projective space.

This image is simple when $D$ is regular, meaning that its colored
points have distinct supports in $C$. For a single
colored point $x=(p,i)$, the intersection
$\overline S_{0,p}\cap\overline T_{-\alpha_i^\vee,p}$ is
$\mathbb P^1$, and $\mathcal H_x$ consists of its two fixed points.
The restriction
isomorphism is simply
\[
\Gamma(\mathbb P^1,\mathcal O(1))
\xrightarrow{\sim}
\Gamma\!\left(
\{0,\infty\},
\mathcal O(1)|_{\{0,\infty\}}
\right)\simeq \b C^2.
\]
For a regular divisor $D=x_1+\cdots+x_n$, factorization and
the restriction isomorphism give
\[
\Gamma\bigl(\mathcal H_D,\mathcal O(1)|_{\mathcal H_D}\bigr)^\vee
\simeq
\bigotimes_{a=1}^n
\Gamma\bigl(Z_{G,x_a},\mathcal O(1)\bigr)^\vee.
\]
Under $Z_{G,D}\simeq\prod_{a=1}^n Z_{G,x_a}\simeq(\mathbb P^1)^n$,
the line bundle $\mathcal O(1)$ corresponds to
$\boxtimes_{a=1}^n\mathcal O_{\mathbb P^1}(1)$.
Its Kodaira embedding is therefore the Segre embedding
\[
Z_{G,D}\simeq(\mathbb P^1)^n
\hookrightarrow
\mathbb P\!\left(
\Gamma\bigl(\mathcal H_D,\mathcal O(1)|_{\mathcal H_D}\bigr)^\vee
\right).
\]

Although the fibers over nonregular divisors are harder to describe
directly, Mirkovi\'c's local projective space construction recovers
the Zastava space as the scheme-theoretic closure of
the regular Segre family in the projective bundle over
$\mathcal H_{C\times I}$ whose fibers the projective spaces
 displayed above for all $D$.

The construction extends beyond groups because the restricted line
bundle $\mathcal O(1)|_{\Gr_T}$, together with its factorization
structure, is determined by the bilinear form $\kappa$ on $X_*(T)$ and can be
constructed directly from $\kappa$. Thus one can start with any
integral symmetric matrix $\kappa$ with even diagonal entries and the
corresponding factorizable line bundle on the colored Hilbert scheme.
The schemes $\mathcal H_D$ and hence the projective bundle used in Kodaira embedding can be defined using $\kappa$ without mentioning $G$.
Taking the closure of the regular Segre family in the projective bundle
defines the generalized Zastava space $Z^{\mathcal H,\kappa}$.
We recall this construction in \S\ref{recollection}.
\begin{rmk}
This construction arose from Mirkovi\'c's work on the
$SL_2$ case of the Kamnitzer--Knutson program on
reconstructing Mirkovi\'c--Vilonen cycles from quivers.
See \cite[\S1.1.6 and \S4]{mirkovic2014loop} for the historical
remarks, and \cite[\S1.6--1.7,12.2]{baumann2021mirkovic}
for developments of the program.
\end{rmk}
\subsubsection{Interval Zastava spaces}

The new construction in this paper extends Mirkovi\'c's Zastava
spaces from effective colored divisors to intervals of arbitrary
colored divisors. Its purpose is to treat all intersections
$\overline S_{-\lambda}\cap\overline T_{-\mu}$ for any
$\lambda\leq\mu$ in $\mathbb Z[I]$ in a uniform construction.

Let $\Div_{C\times I}$ denote the sheaf of $I$-colored relative
Cartier divisors on $C$, and let $\Int$ denote the sheaf of pairs
$(D,E)$ of such divisors with $E-D$ effective, called intervals.
The map $(D,E)\mapsto(D,E-D)$ identifies $\Int$ with
$\Div_{C\times I}\times\mathcal H_{C\times I}$. In fact,
$\Div_{C\times I}$ and $\Int$ are not (ind)schemes,
they are sheaves for the flat topology \cite[3.10.7]{beilinson2004chiral}. Nevertheless,
for each interval $D\leq E$, the divisors $\mathcal I(D,E):=D+\mathcal H_{E-D}\in \Div_{C\times I}$ between $D,E$ define a finite scheme. This generalizes the notion of Hilbert scheme of a finite divisor in the previous section.

To extend the line bundle to arbitrary divisors, we use the Faltings
line bundle $\mathcal L^\kappa$ on $\Bun_T(C)$ \cite{faltings2003algebraic} and its pullback
$\mathbb L^\kappa$ along Abel--Jacobi map
$AJ:\Div_{C\times I}\to\Bun_T(C)$. In the group case,
$\mathcal L^\kappa$ pulls back to
$\mathcal O(1)|_{\Gr_T}$.
The Faltings line bundle depends only on $\kappa$ and does not
require the existence of $G$.
We can define a vector bundle $\b V$ over $\Int$, with fiber
\[
\mathbb V_{(D,E)}
=\Gamma\bigl(\mathcal I(D,E),\mathbb L^\kappa\bigr)^\vee.
\]
Applying the local projective space construction over $\Int$
defines the interval Zastava space $Z=Z^\kappa\subset\mathbb P(\mathbb V)$,
with fibers $Z_{(D,E)}$.

Our main results are as follows.
\begin{thm}
Fix an open embedding $C\hookrightarrow\mathbb A^1_z$. 
\begin{enumerate}
\item(Theorem~\ref{thm:interval-translation})
There are canonical isomorphisms
\[
f_D:\mathbb P\bigl(\mathbb V_{(A,B)}\bigr)
\xrightarrow{\sim}
\mathbb P\bigl(\mathbb V_{(A+D,B+D)}\bigr)
\]
for every interval $A\leq B$ and every colored divisor $D$.
The restrictions give isomorphisms of Zastava spaces: $f_F:Z_{(A,B)}\xrightarrow{\sim}Z_{(A+D,B+D)}$.
\item(Theorem~\ref{thm:zastava-comparison}) 
Let \(Z^{\mathcal H,\kappa}\) denote Mirković’s Zastava space over the colored Hilbert scheme \(\mathcal H_{C\times I}\). Then our interval Zastava space is canonically isomorphic to its pullback along \(\delta:\operatorname{Int}\to\mathcal H_{C\times I}\), where \(\delta(D,E)=E-D\).
\end{enumerate}

\end{thm}

The key to the proof is the identification of the Faltings line bundle
$\mathbb L^\kappa$ with the dual of the line bundle induced by the Heisenberg central extension of \(T_{\mathcal K}\)
(Proposition~\ref{prop:faltings-heisenberg-comparison}).
Multiplication by affine torus group elements realizes translation (i.e. addition) of divisors. 

Finally, inclusions of intervals induce closed embeddings
$Z_{(D,E)}\hookrightarrow Z_{(F,G)}$ when
$(D,E)\le(F,G)$. Thus for any point $p\in C$, we define an indscheme
\[
\mathcal G^\kappa_p
:=\varinjlim_{\lambda\leq\mu}
Z_{(\lambda[p],\mu[p])}.
\]
Allowing the support to vary gives a factorization space over
$\mathcal H_C$. When $\kappa$ is a Cartan matrix of type ADE,
these constructions recover the affine and 
Beilinson--Drinfeld Grassmannians, see \S\ref{subsec:interval-grassmannian}.
\subsubsection{Frenkel--Kac--Segal theorem}

For a simple Lie algebra $\mathfrak g$ of type ADE,
with Cartan matrix $\kappa$, let $L_1(\mathfrak g)$ be the space of basic representations of level 1, and  $L_\kappa$ be the lattice vertex algebra associated to the root lattice.  The Frenkel--Kac--Segal theorem 
\cite{frenkel1980basic,segal1981unitary} states that 
$L_1(\mathfrak g)$ and $ L_\kappa$ are isomorphic as vertex algebras. See the introduction of \cite{zhu2009affine}.

Zhu gave a geometric proof of this isomorphism
\cite[\S3]{zhu2009affine} as follows.
By the affine Borel--Weil theorem, the basic representation $L_1(\f g)$ is realized as the dual space of global sections 
$L_1(\f g)\simeq \Gamma(\Gr_G,\c O_{\Gr_G}(1))^*$. On the other hand, the lattice vertex algebra $L_\kappa$ is
realized as $\Gamma(\Gr_T,O_{\Gr_G}(1)|_{\Gr_T})^*$, see section \ref{Fock}. 
$\Gr_T$ is the \(T\)-fixed point subscheme of $\Gr_G$. Hence it suffices to prove that the restriction of global sections is an isomorphism. Zhu's proof uses Schubert varieties of $\Gr_G$.

As explained above, Mirkovi\'c's Zastava construction yields an
analogue of the Frenkel--Kac--Segal theorem for general $\kappa$
\cite{mirkovic2021loop,mirkovic2023lattice}. In type ADE, it gives
another proof of the classical theorem.

Our interval construction expresses the lattice vertex algebra in
terms of an inductive limit. For any such $\kappa$ and point $p\in C$,
we obtain
\[
\Gamma\bigl(\mathcal G^\kappa_p,\mathcal O(1)\bigr)
\xrightarrow{\sim}
\Gamma\bigl(\Gr_{T,p},\mathbb L^\kappa\bigr)
\simeq L_\kappa^*.
\]
Hence $L_\kappa$ can be realized as a colimit of dual spaces of global sections of Zastava fibers. These isomorphisms are compatible with factorization and give a
geometric realization of $L_\kappa$ through interval Zastava spaces,
see Proposition~\ref{prop:grassmannian-lattice-vertex-algebra}.

\subsection{Related work}
We briefly discuss several works related to the present construction.
The arguments in the main text do not rely on the works discussed below.

\subsubsection{Heisenberg local systems}
Our description of the Faltings line bundle on colored
effective divisors in Section~\ref{subsec:faltings-local-line}
is closely related to the Heisenberg local systems studied in
\cite{bezrukavnikov2006factorizable}, and later \cite{gaitsgory2019metaplectic, gaitsgory2021factorization}.

\subsubsection{BFN Coulomb branches}

The Coulomb branch construction of Braverman--Finkelberg--Nakajima
provides a generalization of affine Grassmannian slices
for any quiver \cite{BFN2}, it is called Coulomb branch of quiver gauge theory with framing. Their compactified Coulomb branches
for unframed quiver gauge theories are identified with Mirkovi\'c's
Zastava by Dong's theorem\cite{dong2026loop}. We do not know how to relate compactified Coulomb branches for framed quiver gauge theories to Mirković’s Zastava spaces. We hope that our construction will help establish such a relationship.

Muthiah and Weekes \cite{muthiah2024fundamental,muthiah2026coulomb} construct inclusions of generalized slices for Kac–Moody quivers under certain hypotheses. These results can be understood as a parallel problem of constructing compatible maps among finite-dimensional pieces in hypothetical affine Grassmannian of Kac--Moody type. If the relationship described in the preceding paragraph can be established, it would be interesting to compare their transition maps with ours.

\subsubsection{Intersection cohomology}

Zastava spaces serve as finite-dimensional models for
semi-infinite flags. Hayash studies the semi-infinite
intersection cohomology sheaf
 on Mirkovi\'c zastava spaces \cite{hayash2025semiinfinite}.

\subsubsection{Other Zastava spaces}

There are versions of Zastava spaces associated
with other curves. Finkelberg--Kuznetsov--Rybnikov study
trigonometric Zastava spaces \cite{finkelberg2018towards}, and 
Finkelberg--Matviichuk--Polishchuk develop elliptic
Zastava spaces \cite{finkelberg2023elliptic}.

\subsubsection{The stack approach}
Bouthier and Vasserot define the affine Grassmannian of a Kac–Moody group using the quotient \([G_{\mathcal K}/G_{\mathcal O}]\), viewed as an étale sheaf in \(\infty\)-stacks. For affine Kac–Moody groups, they establish a geometric Satake equivalence between abelian categories \cite{bouthier2025satake}.
When $\kappa$ is an affine Cartan matrix, it is natural to ask whether the indscheme \(\mathcal G^\kappa\) constructed here admits an embedding into the quotient \([G_{\mathcal K}/G_{\mathcal O}]\).

\subsection{Notation and conventions}
Throughout this paper, for simplicity we work over $\b C$, but the results and proofs extend to general ground rings. 

Let \(I\) be a finite set, and let \(\mathbb Z[I]\) denote the free abelian group with basis \(I\). By a lattice \(L=(\mathbb Z[I],\kappa)\), we mean this group equipped with an even integral symmetric bilinear form \(\kappa\).

Let \(G\) be a connected, simple, simply connected algebraic group. Fix a maximal torus \(T\subset G\) and a Borel subgroup \(B\) containing \(T\), and let \(W\) be the Weyl group. Whenever a group $G$ appears, the lattice $(\b Z[I],\kappa)$ associated to $G$ is the coweight lattice together with the minimal even integral symmetric bilinear form (i.e. normalized bilinear form) on it.

Let $C$ be a smooth connected projective curve. A divisor $D\in\Div_C$ on $C$ is a finite sum $D=\sum_i n_ip_i$, where $n_i\in \b Z$, $p_i\in C$. If $n_i$ is nonnegative for all $i$, then $D$ is said to be \textit{effective}. The moduli space of effective divisors of degree $n$ is the Hilbert scheme of $n$ points of $C$, denoted by $\c H_C^n$. The \textit{Hilbert scheme} $\c H_C$ is the disjoint union $\c H_C=\coprod_{n\ge0} \c H_C^n$. Set $\Div_{C\times I}:=\prod_{i\in I}\Div_C$ (resp. $\c H_{C\times I}:=\prod_{i\in I}\c H_C$), representing $I$-colored (resp. effective) divisors on $C$. Let $\Bun_T(C)$ (resp. {$\Bun_G(C)$}) denote the moduli stacks of \(T\)-bundles and \(G\)-bundles on \(C\), respectively. We omit \(C\) from the notation when no confusion  arises. After choosing an isomorphism $T\simeq \b G_m^I$, the Abel--Jacobi map $AJ: \Div_{C\times I}\to\Bun_T(C)$ sends $D=(D_i)_{i\in I}$ to the direct sum of line bundles $\bigoplus_{i\in I}\c O_C(-D_i)$. We use this negative sign convention throughout.

All closures are scheme-theoretic.
\subsection{Contents}
The paper is organized as follows. Section 2 introduces the Faltings line bundle over $\Bun_T(C)$ and the factorization structure for its Abel--Jacobi pullback to $\Div_{C\times I}$, and computes its restriction to effective colored divisors. Section 3 develops the interval Zastava construction and proves the translation and comparison theorems. Then we  construct generalized affine and Beilinson-–Drinfeld Grassmannians as inductive limits and gives a realization of lattice vertex algebras. Section 4 constructs the Tits group action on generalized affine Grassmannian for simply-laced Kac--Moody Cartan matrices.
\subsection{Acknowledgment}
The author is deeply grateful to his advisor Ivan Mirkovi\'c for teaching him mathematics, continued encouragement, and numerous comments on an early draft of this paper. The author would like to thank Tom Braden for his  comments that improved the writing. The author also thanks Zhijie Dong for useful discussions and pointing out a mistake.
	
\section{Faltings line bundle}
\subsection{Introduction}
This section is about some algebraic geometry preliminaries. The aim is to define the Faltings line bundle \cite{faltings2003algebraic} over $\Bun_T(C)$ for a smooth projective curve $C$, and study its properties.

The definition of Faltings line bundle depends only on the associated lattice $L=(\b Z[I],\kappa)$. In particular, it applies beyond ADE type. We also prove that its pullback  along the Abel--Jacobi map is a factorizable line bundle. 
Hence the Faltings line bundle provides a natural framework to study zastava spaces using techniques in \cite{mirkovic2021loop,mirkovic2023lattice} without using the existence of a group $G$. 
We also compute Faltings line bundle explicitly over the colored Hilbert scheme $\c H_{C\times I}$.

In this section, we first recall the norm line bundle, the $\lambda$-bundle \cite{beilinson2004chiral}, and the Deligne pairing \cite{deligne67determinant}. We then define the Faltings line bundle and study its properties.

\subsection{Norm bundle}\cite{deligne67determinant}\label{normbundle}
Let $D\subset C\times S$ be an effective relative Cartier divisor over test scheme $S$, so that 
$p: D\to S$ is finite and flat. For any line bundle $M$ on $C\times S$, 
the norm line bundle of $M$ along $D$ over $S$ is defined by
\[\Nm_{D/S}(M):=\det p_*(M|_D)\otimes \det\!^{-1} p_*(\c O_D)\in\Pic(S),\]
where $\det$ denotes the determinant line bundle \cite{knudsen1976projectivity}.
On the fiber of a geometric point, if $D=\sum_i n_i p_i,\  p_i\in C$, one has a canonical identification
$\operatorname{Nm}_{D/S}(M)\cong \bigotimes_i M_{p_i}^{n_i}.$ 
For any effective divisors $D_1$ and $D_2$, 
there is a canonical isomorphism
$\operatorname{Nm}_{D_1+D_2/S}(M)
\cong
\operatorname{Nm}_{D_1/S}(M)\otimes \operatorname{Nm}_{D_2/S}(M).$

\subsection{$\lambda$-bundle}
We would like to generalize the notion of norm bundle (for effective divisors) to all divisors. 
\subsubsection{Definition}\label{defn-lambdabundle}
Here we recall the definitions in \cite[3.10.7]{beilinson2004chiral}.
Denote by $\Div_C$ a functor which assigns to each quasi-compact scheme $S$ the group of relative Cartier divisors on $C\times S\to S$ proper over $S$. This is a sheaf with respect to the flat topology.

Let $C$ be a smooth curve and let $p:C\times S\to S$ be the projection. Choose locally on $S$ an effective divisor $D'$ such that $D+D'\ge 0$, and set
\[
\lambda_D
:=
\det\bigl(p_{*}(\mathcal O_{C\times S}(D)/\mathcal O_{C\times S}(-D'))\bigr)
\otimes
\det\!^{-1}\bigl(p_{*}\mathcal O_{D'}\bigr),
\]
where $\mathcal O_{D'}:=\mathcal O_{C\times S}/\mathcal O_{C\times S}(-D')$. 
Both $p_{*}(\mathcal O_{C\times S}(D)/\mathcal O_{C\times S}(-D'))$ and $p_{*}\mathcal O_{D'}$ are vector bundles on $S$, 
and it is shown in \textit{loc. cit.} that $\lambda_D$ is canonically independent 
of the choice of $D'$. Thus one obtains a well-defined (super) line bundle $\lambda$ on 
$\operatorname{Div}_C$. 
For any line bundle $M$ on $C\times S$, the $M$-twisted $\lambda$-bundle is defined by
\[
\lambda_D(M)
:=
\det\bigl(p_{*}((\mathcal O_{C\times S}(D)/\mathcal O_{C\times S}(-D'))\otimes M)\bigr)
\otimes
\det\!^{-1}\bigl(p_{*}(\mathcal O_{D'}\otimes M)\bigr),
\]
which is again independent of $D'$. 

\subsubsection{Relation with norm bundle}
For an effective relative Cartier divisor $D$, the definitions give $\lambda_{D}=\det p_* \c O_D(D)$, $\lambda_D(M)=\det p_*(\c O_D(D)\otimes M)$ by taking $D^\prime=0$. Since $p|_D$ is finite, its pushforward on  quasi-coherent sheaves is exact. Hence $$p_*(\c O_D(D)\otimes_{\c O_D} M)=p_* \left(M|_D \right)\otimes_{p_*\c O_D}p_*\c O_D(D).$$

To compute the determinant, we need the following lemma:
\begin{lem}
Let $A$ be a finite locally free $\c O_S$-algebra, and let $P$ and $Q$ be invertible $A$-modules; then there is a canonical isomorphism
$$\det\!_{\c O_S} \left( P\otimes_{A}Q\right)\otimes\det\!_{\c O_S}^{-1} Q\cong \det\!_{\c O_S}  P\otimes \det\!_{\c O_S}^{-1}A$$
\end{lem}
\begin{proof}
Let $q$ be a generator of $Q$ as $A$-module. Consider the exact sequence
$$0\to P\xrightarrow{f} P\otimes_AQ\oplus A\xrightarrow{g}Q\to 0,$$
$f(p)=(p\otimes q,0), p\in P,\quad  g(-,a)=aq$. Taking determinants give $\det(P)\otimes \det(Q)\simeq \det(P\otimes_AQ\oplus A)\simeq \det(P\otimes_AQ)\otimes \det(A)$.
To check the isomorphism is canonical, take $uq$ as another generator, where $u$ is a unit of $A$. Then both sides in the formula multiply by $\det(u:A\to A)$, which preserves the formula. 
\end{proof}
Take $P=p_* \left(M|_D \right), Q=p_*\c O_D(D), A=p_*\c O_D$, then	 by the lemma above 
$$\lambda_D(M)\otimes \lambda_D^{-1}=\det \left( P\otimes_{A}Q\right)\otimes\det\!^{-1} Q\cong \det  P\otimes \det\!^{-1}A=\op{Nm}_D(M).$$

\subsection{Deligne pairing \cite{deligne67determinant}}
Let $C$ be a smooth projective curve, and $S, p:C\times S\to S$ be as above. For line bundles $\c L$ and $\c M$ on $C\times S$, define their \emph{Deligne pairing} as the following line bundle on $S$:
\[\langle \c L,\c M\rangle
:={}\det Rp_*(\c L\otimes \c M)
\otimes\det Rp_*\c O_{C\times S}
\otimes\det\!^{-1}Rp_*\c L
\otimes\det\!^{-1}Rp_*\c M.
\]
The Deligne pairing is symmetric, bimultiplicative (meaning $\langle \c L,\c M\rangle\otimes\langle \c L^\prime, \c M\rangle\cong \langle \c L\otimes \c L^\prime,\c M\rangle$), and compatible with base change. In particular for any $\c L$,
$\langle \c L,\c O_{C\times S}\rangle\cong\c O_S.$

The next proposition identifies the Deligne pairing with the 
$\lambda$-bundle introduced above.

\begin{prop}
	\label{prop:deligne-lambda}
	Let $D\in\Div_C(S)$ be a relative Cartier divisor (not necessarily effective), and  $\c L$ be a line bundle on $C\times S$. Then there is a canonical isomorphism
	\[	\langle\c O_{C\times S}(D),\c L\rangle
	\cong	\lambda_D(\c L)\otimes\lambda_D^{-1}.	\]
 As a corollary, if $D$ is effective, we have
	$\langle\c O_{C\times S}(D),\c L\rangle
	\cong \lambda_D(\c L)\otimes\lambda_D^{-1}\cong 
	\Nm_{D/S}(\c L)$.
\end{prop}

\begin{proof}
	Write  $\c O_{C\times S}$ simply as $\c O$. Choose an effective relative Cartier divisor $D'$ such that $D+D'$ is effective as in the definition of $\lambda$-bundle in \ref{defn-lambdabundle}. 
	Consider short exact sequences
	\[
	0\to\c O(-D')\to\c O(D)\to \c O(D)/\c O(-D')\to0,
	\qquad
	0\to\c O(-D')\to\c O\to \c O_{D'}\to0.
	\]
	Taking determinants we get
	$$\det Rp_* \c O(-D')\otimes\det Rp_*\left(\c O(D)/\c O(-D') \right) \cong \det Rp_*\c O(D),\qquad
	\det Rp_*\c O(-D')\otimes \det Rp_*\c O_{D'}\cong \det Rp_*\c O.$$
	Hence $$\det Rp_*\c O(D)\otimes \det\!^{-1} Rp_*\c O\cong \det Rp_* (\c O(D)/\c O(-D^\prime))\otimes\det\!^{-1}Rp_*\c O_{D^\prime}=\lambda_{D}  .$$
	Similarly by considering the exact sequence tensoring with $\c L$
	we get $$\det Rp_*(\c L\otimes\c O(D))\otimes \det\!^{-1} Rp_*(\c L\otimes\c O)\cong \det Rp_* (\c L\otimes\left(\c O(D)/\c O(-D^\prime) \right) )\otimes\det\!^{-1}Rp_*(\c L\otimes\c O_{D^\prime})=\lambda_{D}(\c L) .$$
	Combining the second formula with the inverse of the first we get the stated proposition.
\end{proof}

\subsection{Faltings line bundle \cite{faltings2003algebraic}}
\subsubsection{Definition}
Let $C$ be a smooth connected projective curve, and $T$ be a torus. We can define the tautological $T$-bundle $\c T$ over $\Bun_T(C)\times C$, whose fiber $\c T_{(P,x)}$ over pair $(P,x)\in \Bun_T(C)\times C$ is $\c T_{(P,x)}=P_x$. 
Let $\{e_i\}_{i\in I}\subset X_*(T)=\Hom(\b{G}_m,T)$ be a basis of coweight lattice, $\{f_i\}_{i\in I}\subset X^*(T)=\Hom(T,\b{G}_m)$ be the dual basis. Hence $\{f_i\}$ identifies $T$ with $\b G_m^n$.

Let $\kappa=(\kappa_{ij})$ be the symmetric bilinear form on $X_*(T)$ defined by $\kappa_{ij}=\kappa(e_i, e_j)$. For instance, if $\{e_i\}$ are simple coroots in ADE case, then $\kappa$ is the corresponding Cartan matrix.

Define $\c{L}_i:=\c T\stackrel{T}{\times}\b C_{f_i}$ to be the line bundle over $\Bun_T(C)\times C$ given by the $i$-th component of $\c T$, where $\b C_{f_i}$ is the 1-dimensional representation of $T$ given by $f_i$.
The \textit{Faltings line bundle} $\c{L}^\kappa$ over $\Bun_T(C)$ is defined as follows in terms of Deligne pairings: 
\[
\mathcal{L}^\kappa:=\bigotimes_{i\le j}\langle \mathcal{L}_i,\mathcal{L}_j\rangle^{-\kappa_{ij}^\prime},
\]
here the matrix $(\kappa^\prime_{ij})$ is obtained from $\kappa$ by halving its diagonal entries; in other words, $\kappa_{ij}^\prime=\frac{1}{1+\delta_{ij}}\kappa_{ij}$. Then $\c L^\kappa$ is independent of the choice of basis, which can be checked directly using bimultiplicative property of the Deligne pairing.

\subsubsection{The case of a group $G$}
 Let $G$ be a simple and simply connected group, with $T$  a maximal torus and $\kappa$ the minimal (i.e. normalized) even integral bilinear form of $G$, $\c O_{\Gr_G}(1)$ be the rigidified ample line bundle generator.
 For any point $p\in C$, let 
 $\c O_p:=\widehat{\c O}_{C,p}$, $\c K_p:=\operatorname{Frac}(\c O_p)$,  $\Gr_{G,p}:=G(\c K_p)/G(\c O_p)\simeq\Gr_G$, $\Gr_{T,p}:=T(\c K_p)/T(\c O_p)$.
 Let  $\pi_G: \Gr_{G,p}\to \Bun_G(C)$ and  $\pi_T: \Gr_{T,p}\to \Bun_T(C)$ be the forgetful maps, then
 the line bundle  $\pi_T^*\c L^\kappa$ is the restriction of  $\c O_{\Gr_G}(1)$ to $\Gr_T$. See \cite[Theorem 17]{faltings2003algebraic}.
 
  Faltings also proves $\c L^\kappa$ is the pullback of a canonically defined line bundle $\c L_G$ over $\Bun_G(C)$ to $\Bun_T(C)$. Here $\c L_G$ is the rigidified ample line bundle generator over $\Bun_G(C)$, and pullback $\pi_G^*\c L_G$ is the rigidified ample generator $\c O_{\Gr_G}(1)$.

However, the construction of $\c{L}^\kappa$ itself relies only on the matrix $(\kappa_{ij})$ and does not require the existence of a corresponding group $G$.

\subsubsection{A representation-theoretic interpretation of global sections}
\label{Fock}
We record an interpretation of the dual space of global sections in terms of Weyl modules. See \cite[\S17.1]{GaitsgoryCFT} for reference.

Let $\wh{\f t}_\kappa=\f t\otimes\c K\oplus\b Cc$ be the Heisenberg algebra associated with $\kappa$. For $\mu\in\f t^*$, let $\b C_\mu$ denote the one-dimensional representation of $\f t\otimes\c O\oplus\b Cc$ on which $\f t\otimes\c O$ acts through $\f t$ by $\mu$, and $c$ acts by 1. Define the Weyl module
\[
V^\mu
:=
\Ind^{\wh{\f t}_\kappa}_{\f t\otimes\c O\oplus\b Cc}
\b C_\mu.
\]
For any $\kappa$ including degenerate $\kappa$, there is an isomorphism of $\wh{\f t}_\kappa$-modules
\[
\Gamma(\Gr_T,\pi_T^*\c L^\kappa)^*
\simeq
\bigoplus_{\lambda\in X_*(T)}
V^{\kappa(\lambda,-)},
\]
here we view $\kappa(\lambda,-)$ as an element in $\f t^*$.

The proof of \cite[Proposition 2.1.1]{zhu2009affine} uses the irreducibility of the corresponding Fock modules, which relies on nondegeneracy of the bilinear form $\kappa$. In the current paper we allow $\kappa$ to be degenerate.

\subsubsection{Comparison with determinant line bundle}
For $G=\SL_n$, we verify fiberwise the comparison between
$\c O_{\Gr_G}(1)|_{\Gr_T}$ and the pullback of the Faltings
line bundle, using the determinant description of
$\c O_{\Gr_G}(1)$. For general $G$, see
\cite[Remark 2.4.3]{zhu2016introduction}.
This calculation is not used later, so readers may feel free to skip
this subsection.

Let $G=\SL_n$, let $T\subset G$ be the diagonal maximal torus, and let 	$\lambda=(\lambda_1,\ldots,\lambda_n)\in X_*(T)$, so that	$\sum_i\lambda_i=0$.
	Choose a smooth projective curve $C$, a point $p\in C$, and a local 	coordinate $z$ at $p$ such that $z(p)=0$. Write 	$\c O_p:=\widehat{\c O}_{C,p}$, $\c K_p:=\operatorname{Frac}(\c O_p)$, and $\Gr_{G,p}:=G(\c K_p)/G(\c O_p)\simeq\Gr_G$.
	
	Let $L_0=\c O_p^n$, $L_\lambda=z^\lambda L_0
	=\bigoplus_{i=1}^n z^{\lambda_i}\c O_pe_i.$
	Recall $\c O(1)$ is the inverse of the 
	determinant line bundle. Its fiber at $z^\lambda$ is 
	\[
	\begin{aligned}
		\left.\c O(1)\right|_{z^\lambda}
		&=
		\det\!^{-1}\bigl(L_\lambda/(L_\lambda\cap L_0)\bigr)
		\otimes
		\det\bigl(L_0/(L_\lambda\cap L_0)\bigr) \\
		&\cong
		\bigotimes_{\lambda_i<0}
		\det\!^{-1}\left(\frac{z^{\lambda_i}\c O_pe_i}{\c O_pe_i}\right)
		\otimes
		\bigotimes_{\lambda_i>0}
		\det\left(\frac{\c O_pe_i}{z^{\lambda_i}\c O_pe_i}\right).
	\end{aligned}
	\]
	
	Let $\f m_p\subset\c O_p$ be the maximal ideal. The local coordinate
	$z$ is an element of $\f m_p$, and its image in
	$\f m_p/\f m_p^2= T_p^*C$ is a  cotangent
	vector. More generally, $z^j\c O_p/z^{j+1}\c O_p\cong (T_p^*C)^{\otimes j}$
	for every $j\in\b Z$. Applying this and using the multiplicativity  of determinants in short exact sequence we get
	\[
	\left.\c O(1)\right|_{z^\lambda}
	\simeq
	(T_p^*C)^{\otimes\frac{1}{2}\kappa(\lambda, \lambda)},
	\qquad
	\kappa(\lambda, \lambda)
	=
	\sum_{i=1}^n(\lambda_i^2-\lambda_i)
	=\sum_{i=1}^n\lambda_i^2.
	\]
	
	We now recover the same description from the Faltings line bundle.
	The image of colored divisor $D=(\lambda_1p,\cdots, \lambda_np)$ under Abel--Jacobi map corresponds to $z^\lambda$. 
	Take the usual simple coroots $e_i=(0,\cdots, 1,-1,\cdots,0)\in X_*(T)$ as basis, under this basis $\lambda$ has coordinates $(m_1,\cdots, m_{n-1})$, where $m_i=\sum_{j=1}^{i}\lambda_j$.
	The Faltings line bundle fiber over $D$ is
	\begin{align*}
		\b L^\kappa_D=&\left(\bigotimes_{i=1}^{n-1}\left\langle \c L_i,\c L_i\right\rangle |_D\right)^{\otimes -1}
		\otimes \left(\bigotimes_{i=1}^{n-2}\left\langle \c L_i,\c L_{i+1}\right\rangle|_D \right)\\
		=&\left(\bigotimes_{i=1}^{n-1}\left\langle \c O(-m_ip),\c O(-m_ip)\right\rangle \right)^{\otimes -1}
		\otimes \left(\bigotimes_{i=1}^{n-2}\left\langle \c O(-m_ip),\c O(-m_{i+1}p)\right\rangle \right)\\
		=&\left(\bigotimes_{i=1}^{n-1} (T_p^*C)^{\otimes m_i^2} \right)
		\otimes \left(\bigotimes_{i=1}^{n-2}(T_p^*C)^{\otimes -m_im_{i+1}}  \right) &\text{(Bimultiplicativity)}\\
		=& (T_p^*C)^{\otimes \left(\sum_{i=1}^{n-1}m_i^2- \sum_{i=1}^{n-2}m_im_{i+1}\right) }
	\end{align*}
	One can check directly that $\sum_{i=1}^{n-1}m_i^2- \sum_{i=1}^{n-2}m_im_{i+1}=\frac{1}{2}\sum_{i=1}^n \lambda_i^2$, hence the two descriptions agree.

\subsubsection{Invariance of $\c L^\kappa$ under inversion}
We will use the following lemma later.
\begin{lem}\label{dual}
	Let $\iota:\Bun_T(C)\longrightarrow\Bun_T(C)$
	be the involution induced by the inversion 
	\(t\mapsto t^{-1}\) of \(T\). Then there is a canonical isomorphism $\iota^*\c L^\kappa\xrightarrow{\sim}\c L^\kappa$.

\end{lem}

\begin{proof}
	Recall $\c L_i:=\c T\mathop{\times}\limits^T\mathbb C_{f_i}$ is the corresponding line bundle on \(\Bun_T(C)\times C\) used in $\c L^\kappa$, where $\b C_{f_i}$ is the 1-dimensional representation of $T$ by $f_i$.	Since \(\iota\) is induced by inversion on \(T\), there is a	canonical isomorphism
	\[(\iota\times\id_C)^*\c L_i\cong \c T\mathop{\times}\limits^T\mathbb C_{-f_i}= \c L_i^*.\]
	
	For any two line bundles \(\c L_i,\c L_j\), the
	bimultiplicative property of the Deligne pairing gives a canonical 	isomorphism
	$\langle\c L_i^*,\c L_j^*\rangle	\cong
\langle\c L_i,\c L_j^*\rangle^*\cong
\langle\c L_i,\c L_j\rangle.$ Hence $\iota^*\c L^\kappa\cong \bigotimes_{ i\le j}\langle\c L_i^*,\c L_j^*\rangle^{\otimes-\kappa'_{ij}}\cong \c L^\kappa$.
\end{proof}

\subsection{The Faltings line bundle as factorizable line bundle}
\subsubsection{Factorization structure}
Let $M$ be a line bundle over $\Div_{C\times I}$ (resp. $\c H_{C\times I}$). A \emph{factorization structure} on $M$
contains the following data \cite[3.10.7]{beilinson2004chiral}:
\begin{enumerate}
\item Factorization: for any pair of disjoint colored Cartier divisors (resp. disjoint colored effective Cartier divisors)  $D,E$, there is a canonical isomorphism $c_{D,E}: M_{D}\otimes M_{E} \xrightarrow{\sim}  M_{D+E}$.
 Here disjointness of colored divisors $D=\{D_i\}_{i\in I}$ and $E=\{E_i\}_{i\in I}$ means $\supp D_i\cap \supp E_j=\emptyset$ for any $i,j$.
\item Consistency: $c$ is compatible with base change, associativity, commutativity. 
\end{enumerate}
If  the condition (1) is strengthened by requiring such a family of isomorphisms $c_{D,E}$ to exist for any $D,E$ not necessarily disjoint, we say $M$ carries a \textit{strong factorization structure}.

In \textit{loc. cit.}, (strong) factorizable line bundles are classified by \textit{theta data}. The classification is also discussed in \cite{tao2021central}.
\begin{ex}
As we have seen in section \ref{normbundle}, the norm bundle over the colored Hilbert scheme $\c H_{C\times I}$ carries a natural strong factorization structure. It is proved in \cite[3.10.7]{beilinson2004chiral} that $\lambda$-bundle carries a natural factorization structure as a super line bundle on $\Div_{C\times I}$.
\end{ex}
\subsubsection{Factorization structure on the Faltings line bundle}\label{faltingslocal}
We now formulate the factorization structure of the Faltings line bundle in terms of divisor presentations.

Recall the Abel--Jacobi map
$AJ:\Div_{C\times I}\to \Bun_T$ sending $
D=\{D_i\}$ to $\bigoplus_i\mathcal O_C(-D_i).$
Denote the pullback of Faltings line bundle $\c L^\kappa$ to $\Div_{C\times I}$ by $$\b L^\kappa := AJ^*\mathcal L^\kappa.$$
From now on we work with the Faltings line bundle $\b L^\kappa$ on $\Div_{C\times I}$ rather than $\c L^\kappa$, and we sometimes suppress $\kappa$ from the notation when no confusion arises.

We begin with the basic triviality property of the Deligne pairing.
\begin{lem}\label{disjoint-trivial}
	Let $D,E\in \Div_C(S)$ be disjoint relative Cartier divisors. There is a canonical isomorphism
	$\langle \mathcal O(D),\mathcal O(E)\rangle \cong \mathcal O_S.$
\end{lem}

\begin{proof}
	Recall the results in Proposition  \ref{prop:deligne-lambda},
	if $D$ is effective, then $\langle \mathcal O(D),\mathcal O(E)\rangle \cong \operatorname{Nm}_D(\mathcal O(E)).$
	Since $D$ and $E$ are disjoint, the restriction $\mathcal O(E)|_D$ is canonically trivial, hence the norm is canonically trivial.
	
	For a general divisor $D$, from the proposition
$	\langle \mathcal O(D),\mathcal O(E)\rangle \cong \lambda_D(\c O(E))\otimes \lambda_D^{-1},$
	holding for any effective divisor $D^\prime$ such that $D+D^\prime\ge 0$ in the definition of $\lambda_{D}$.
	Hence we can choose $D'$ to be supported on $D$, then $\lambda_{D}(\c O(E))$ depends only on the restriction of $\mathcal O(E)$ to the support of $D$, hence $\lambda_{D}(\c O(E))$ is canonically isomorphic to $\lambda_{D}$.
\end{proof}
Now we formulate the factorization structure of $\b L^\kappa$.
\begin{prop}[Factorization]\label{locality}
There is a natural factorization structure on $\b L^\kappa$. 
\end{prop}

\begin{proof}
For disjoint colored divisors $D=\{D_i\},E=\{E_i\}\in \Div_{C\times I}$, we need to construct a system of canonical isomorphisms
\[c_{D,E}:\quad	 \b L_D\otimes \b L_E\xrightarrow{\sim}\b L_{D+E}.\]
		
By definition $\b L_i|_D:=AJ^*\c L_i=\c O(-D_i)$, hence
	\[	\b L_D	=	\bigotimes_{i\le j}	\langle \mathcal O(-D_i),\mathcal O(-D_j)\rangle^{\otimes -\kappa'_{ij}},
	\;	\b L_E	=\bigotimes_{i\le j}
	\langle \mathcal O(-E_i),\mathcal O(-E_j)\rangle^{\otimes -\kappa'_{ij}},\]
	\[\b L_{D+E}
	=
	\bigotimes_{i\le j}
	\langle \mathcal O(-D_i-E_i),\mathcal O(-D_j-E_j)\rangle^{\otimes -\kappa'_{ij}}.
	\]
	
	By bimultiplicativity of Deligne pairing,
	\[		\langle \mathcal O(-D_i-E_i),\mathcal O(-D_j-E_j)\rangle
	\cong
	\langle \mathcal O(-D_i),\mathcal O(-D_j)\rangle
	\otimes
	\langle \mathcal O(-D_i),\mathcal O(-E_j)\rangle 
	\otimes
	\langle \mathcal O(-E_i),\mathcal O(-D_j)\rangle 
	\otimes
	\langle \mathcal O(-E_i),\mathcal O(-E_j)\rangle.
	\]
	
	Since $D$ and $E$ are disjoint, all mixed terms
	$	\langle \mathcal O(-D_i),\mathcal O(-E_j)\rangle,
	\langle \mathcal O(-E_i),\mathcal O(-D_j)\rangle$
	are canonically trivial by Lemma~\ref{disjoint-trivial}. Therefore
$\langle \mathcal O(-D_i-E_i),\mathcal O(-D_j-E_j)\rangle
\cong
\langle \mathcal O(-D_i),\mathcal O(-D_j)\rangle
\otimes
\langle \mathcal O(-E_i),\mathcal O(-E_j)\rangle.$
\end{proof}

\subsection{The Faltings line bundle on effective divisors}
\label{subsec:faltings-local-line}

In this subsection we compute the Faltings line bundle $\b L^\kappa$ explicitly in terms of the factorizable line bundle \(\c O(-\sum_{i\le j}\kappa_{ij}^\prime\Delta_{ij})\) over colored Hilbert scheme $\c H_{C\times I}$. An equivalent formulation has already been announced in \cite[lemma~17.3]{GaitsgoryCFT} in terms of Heisenberg central extensions. Since we do not know a written proof of this identification, we include one here for completeness.

\begin{prop}\label{pairing-delta}
		Let $\c H^n_C\simeq \Sym^n(C)$ be the Hilbert scheme of $n$ points on a smooth curve $C$, and let \(\Delta_n\subset\c H^n_C\) be the discriminant divisor. Also write $\mathcal D_n:=\{c\in C,D\in \c H_C^n\mid c\in D\}\subset C\times \c H^n_C$ for the universal family, and let $q: C\times \c H^n_C\to C$ be the projection. Then there is a canonical isomorphism in $\Pic(\c H^n_C)$:
	\begin{equation}
		\label{eq:self-pairing-universal-divisor}
		\bigl\langle
		\mathcal O(\mathcal D_n),\mathcal O(\mathcal D_n)
		\bigr\rangle
		\simeq
		\operatorname{Nm}_{\mathcal D_n/\c H^n_C}(q^*T_C)
		\otimes\mathcal O_{\c H^n_C}(\Delta_n).
	\end{equation}
\end{prop}
\begin{proof}
Let $\pi: C^n\to \c H^n_C=C^n/S_n$ be the quotient map sending $(x_1,\ldots,x_n)$ to $x_1+\cdots+x_n$. Let \(p_i:C^n\to C\) be the \(i\)-th projection and let
\(\Gamma_i\subset C\times C^n\) be the divisor $\{(x_0,x_1,\cdots, x_n)\in C^{n+1}\mid x_0=x_i\}$. The pullback divisor $\pi^*\c D_n$ is then
$\Gamma_1+\cdots+\Gamma_n$. By bimultiplicativity of Deligne pairing, we have
$$\pi^*\left\langle \c O(\c D_n), \c O(\c D_n)\right\rangle\cong \bigotimes_{i,j=1}^n\left\langle \c O(\Gamma_i), \c O(\Gamma_j)\right\rangle.$$
We calculate the right side in the following steps:
\begin{enumerate}
\item The case $i=j$. 
Consider the map $\Gamma_i \overset{f}\hookrightarrow C^{n+1}\xrightarrow{\id_C\times p_i}C^2$. Let $\Delta\subset C^2$ be the diagonal, then 
$$\c O(\Gamma_i)|_{\Gamma_i}=\left( (\id_C\times p_i)\circ f \right) ^*\c O_{C^2}(\Delta)=p_i^*(\c O_{C^2}(\Delta)|_{\Delta})\cong  p_i^*T_C .$$
Because
\(\Gamma_i\to C^n\) is an isomorphism, 
$\bigl\langle\mathcal O(\Gamma_i), \mathcal O(\Gamma_i)\bigr\rangle
=\op{Nm}_{\Gamma_i}(\c O(\Gamma_i))
\cong\c O(\Gamma_i)|_{\Gamma_i}= p_i^*T_C.$
\item The case $i\neq j$. After identification $\Gamma_i\cong C^n$, the restriction $\c O(\Gamma_j)|_{\Gamma_i}$ is isomorphic to $\c O_{C^n}(\Gamma_{ij})$, where $\Gamma_{ij}\subset C^n$ is given by the equation that $i$-th and $j$-th coordinates are equal.
\item Note that $\pi^*\Delta_n=2\sum_{i<j}\Gamma_{ij}$. This is because the diagonal divisor in Hilbert scheme is given by quadratic equations in Cartesian coordinates, which can be checked locally in the case $n=2$ as follows. 

When $n=2$, $C=\b A^1$, $x_1,x_2$ are the coordinates of $C^2$, with diagonal divisor $\Gamma_{12}=(x_1-x_2=0)$. On the other hand, $\Sym^2(C)=\b A^2_{x_1+x_2,x_1x_2}$, with  $\pi^*\Delta_2=\left((x_1+x_2)^2-4x_1x_2=0 \right)=2\Gamma_{12} $ .

\end{enumerate}
Combining the three points above, we get
\begin{align*}
\bigotimes_{i,j=1}\left\langle \c O(\Gamma_i), \c O(\Gamma_j)\right\rangle=&
\bigotimes_{i}\left\langle \c O(\Gamma_i), \c O(\Gamma_i)\right\rangle
\otimes
\left(\bigotimes_{i<j}\left\langle \c O(\Gamma_i), \c O(\Gamma_j)\right\rangle \right)^{\otimes2} \\
=&\bigotimes_{i} p_i^*T_C\otimes \c O_{C^n}\left( 2\sum_{i<j}\Gamma_{ij}\right) \\
=&\bigotimes_{i} p_i^*T_C\otimes \c O_{C^n}\left( \pi^*\Delta_n\right) 
\end{align*}
The isomorphism constructed above is natural and $S_n$-equivariant. Hence it descends along $\pi$, proving the proposition.
\end{proof}

Now we turn to the mixed color situation.
Let $\c D_i\subset C\times \c H_{C\times I}$ be the universal family recording the $i$-th color, $q: C\times \c H_{C\times I}\to C$ be the projection. Write $\Delta_{ij}\subset \c H_{C\times I}$ for the discriminant divisor and $\Delta_{ij}$ the incidence divisor for $i\neq j$ in $\c H_{C\times I}$. From the proposition above, we know that
$$\left\langle \c O(\c D_i),\c O(\c D_i)\right\rangle
\cong \operatorname{Nm}_{\mathcal D_i/\c H_{C\times I}}(q^*T_C)
\otimes\mathcal O_{\c H_{C\times I}}(\Delta_{ii}).$$
Similarly, for $i\neq j$ we have
$\left\langle \c O(\c D_i),\c O(\c D_j)\right\rangle
\cong \mathcal O_{\c H_{C\times I}}(\Delta_{ij}).$ Here the self-intersection contribution from case (1) in the proof above is absent.
Recall $\b L^\kappa= AJ^*\bigotimes_{1\le i\le j\le n}\langle \mathcal{L}_i,\mathcal{L}_j\rangle^{-\kappa_{ij}^\prime}$, the restriction of $AJ^*\c L_i$ over $\c H_{C\times I}$ is just $\c O_{C\times \c H_{C\times I}}(-\c D_i)$. 

From all the discussions above in this subsection, we conclude the following theorem:
\begin{thm}[Faltings line bundle on effective divisors]\label{thm:faltings-on-colored-hilbert}
With the notation above, over \(\mathcal H_{C\times I}\) there is a canonical isomorphism
	$$\left.\mathbb L^\kappa\right|_{\mathcal H_{C\times I}}	\cong		\left(		\bigotimes_{i\in I}		\operatorname{Nm}_{\mathcal D_i/\mathcal H_{C\times I}}
	(q^*\omega_C)^{\otimes\kappa^\prime_{ii}}\right)
	\otimes	\mathcal O_{\mathcal H_{C\times I}}(-\sum_{i\le j}\kappa^\prime_{ij}\Delta_{ij}).\qquad\qed$$
	
\end{thm}
\begin{rmk}\label{thetadata}
	$\theta$-data introduced in \cite[3.10.3]{beilinson2004chiral} determines a factorization line bundle up to isomorphism.,For the Faltings line bundle $\b L=\b L^\kappa$, its $\theta$-datum is the triple $\left(
	-\kappa,\,
	\{\lambda^\gamma\}_{\gamma\in L},\,
	\{c^{\gamma,\delta}\}_{\gamma,\delta\in L}
	\right)$ given as follows.
	Put	$j_\gamma:C\to\Div_{C\times I}$, $x\mapsto\gamma[x]$, we have
	$\lambda^\gamma:=j_\gamma^*\b L
	\simeq\omega_C^{\otimes \frac12\kappa(\gamma,\gamma)}$. The multiplication maps are
	\[
	c^{\gamma,\delta}:
	\lambda^\gamma\otimes\lambda^\delta
	\xrightarrow{\sim}
	\lambda^{\gamma+\delta}\otimes
	\omega_C^{\otimes-\kappa(\gamma,\delta)}.
	\]
	Let $\Delta\subset C^2$ be the diagonal and set
	$j_{\gamma,\delta}(x,y):=\gamma[x]+\delta[y]$.
	The map $c^{\gamma,\delta}$ is obtained by restricting the corrected
	factorization isomorphism
	\[
	\lambda^\gamma\boxtimes\lambda^\delta
	\xrightarrow{\sim}
	j_{\gamma,\delta}^*\b L\otimes
	\mathcal O_{C^2}\bigl(\kappa(\gamma,\delta)\Delta\bigr)
	\]
	to $\Delta$, using
	$\mathcal O_{C^2}(\Delta)|_\Delta\simeq\omega_C^{-1}$.
	One can check directly that these maps are associative and satisfy
	\[
	c^{\gamma,\delta}
	=
	(-1)^{\kappa(\gamma,\delta)}
	c^{\delta,\gamma}\circ\sigma,
	\]
	where $\sigma$ exchanges the two tensor factors.

\end{rmk}
\section{Zastava for intervals}
\subsection{Recollection of the Mirkovi\'c Zastava construction}\label{recollection}

The main reference for this subsection is \cite[section 2]{mirkovic2021loop}.

\subsubsection{Hilbert scheme of divisors}
\label{Hilbertdivisor}
For a colored effective divisor $D=(D_i)_{i\in I}\in \c H_{C\times I}$, the 
Hilbert scheme $\c H_D$  is the Hilbert scheme of closed subschemes of the finite scheme $D$, i.e. consisting of all 
effective colored divisors $E$ with $0\le E\le D$ with multiplicities.
\footnote{It is also called Grassmannian $\Gr(D)$ of effective divisor $D$ in \cite{mirkovic2021loop,mirkovic2023lattice,dong2026loop}.} It is in general not reduced.
In other words, the scheme structure is given by Hilbert scheme of points in $D$:
$$\c H_D=\prod_{i\in I} \c H_{D_i}=\prod_{i\in I} \coprod_{j=0}^{\deg(D_i)}\c H^j_{D_i}.$$ The length of 
$\c H_{D_i}$ is $2^{\deg(D_i)}$, the length of $\c H_D$ is $2^{\deg(D)}:=2^{\sum \deg(D_i)}$.
For example, when $D=n[p]$, $(\c H_D)_{\text{red}}=\{0,[p],\cdots,n[p]\}$. For $i\le n$, the multiplicity of $i[p]$ is $\binom{n}{i}$.  

\subsubsection{Local vector bundle $\c V^\kappa$}
\label{c V}
Denote by $\Gr(\c D/\c H)$ the incidence subscheme of $\c H_{C\times I}\times \c H_{C\times I}$ parameterizing pairs of effective divisors $(D',D)$ satisfying $D'\le D$.
Consider the correspondence $$\c{H}_C\xleftarrow{q} \Gr(\c D/\c H)\xrightarrow{p} \c{H}_C$$ given by $D^\prime\mapsfrom (D^\prime\le D)\mapsto D$. Recall the Faltings line bundle $\mathbb L=AJ^*\c L^\kappa$ on $\Div_{C\times I}$ in section \ref{faltingslocal}, define sheaf 
\footnote{Here we use the Faltings line bundle rather than $\c O(-\sum_{i\le j}\kappa'_{ij}\Delta_{ij})$ over $\c H_{C\times I}$. This differs from \cite{mirkovic2021loop}, they are isomorphic after choosing a local coordinate in formal disc. See section \ref{subsec:faltings-local-line} for comparison details.}
$$\c V^\kappa=(p_*q^*\b{L}^\kappa)^*.$$

Since $p$ is finite locally free (see \textit{loc. cit.}, section 2.4.2), $\c V^\kappa$ is a vector bundle. Using flat base change, one can check the fiber $\c{V}_D$ over $D\in\c{H}_C$ is $\c{V}_D=\Gamma(\c H_D,\b{L}^\kappa)^*$. 

The terminology ``local" will be explained later in section \ref{local}.
\subsubsection{Mirkovi\'c Zastava space}
The construction is given by the following steps:
\begin{enumerate}
	\item
	For a point $p=(a,i)\in C\times I$,
	$\c{V}_p=\Gamma(\c H_p,\b{L}^\kappa)^*=\b{L}_0^*\oplus \b{L}_p^*$ is a 2-dimensional vector space. Define the fiber of \textit{Zastava space} $Z^{\c H}_p$ to be the projectivization $Z^{\c H}_p:=\b{P}(\c{V}_p)$ which is isomorphic to a projective line.
	
	\item	
	Now we define the Zastava space 
	 $Z^{\c H}:=\b{P}^{\text{loc}}(\c{V})$ over $\c{H}_{C\times I}$. 
	 \footnote{The superscript $\c H$ is to emphasize that the base is Hilbert scheme $\c H_{C\times I}$, and to distinguish the interval version introduced later.} 
	 Its fiber over a regular divisor $D\in \c{H}_{C\times I}^{\text{reg}}=\c{H}_{C\times I}\sm \bigcup_{i\le j}\Delta_{ij}$ is defined as the image under the Segre embedding $$Z^{\c H}_D=\im\left( \prod_{p\in D}\b{P}(\c{V}_p)\hookrightarrow\b{P}\left( \bigotimes \c{V}_p\right) \right) .$$ 
	\item 
	The whole space $Z^{\c H}=\b{P}^{\text{loc}}(\c V)$ is defined as the closure of its regular part $(Z^{\c H})^{\text{reg}}$ in $\b{P}(\c{V})$. Note this construction $\b{P}^{\text{loc}}(\c{V})$ makes sense for any local vector bundle $\c{V}$ with the same procedure. We call  $\b{P}^{\text{loc}}(\c{V})$ the \textit{local projective space} of $\c{V}$.
\end{enumerate}
\begin{rmk}
	The Faltings line bundle $\b L^\kappa$ depends only on the
	lattice $(\b Z[I],\kappa)$ and is independent of the choice of basis.
	However, the Zastava construction $Z^{\c H,\kappa}$ uses the chosen basis to
	specify effective colored divisors, hence it depends on
	$\kappa$ as a \textit{symmetric matrix} in a chosen basis, rather than only as a bilinear form. In particular, integrally congruent
	symmetric matrices need not yield isomorphic Zastava spaces,
	see Remark~\ref{rmk:affine-a-congruence}.
\end{rmk}
We record some examples of two-point Zastava fibers, referring to
\cite{mirkovicnoteszastava} for the calculations.

\begin{ex}
	\label{ex:zastava-two-point-fibers}
	Let $p\in C$.
	
	For distinct colors $i\ne j$, there are homogeneous coordinates
	$[X_0:X_i:X_j:X_{ij}]$ on $\b P^3$ in which
	\[
	Z^{\c H}_{(p,i)+(p,j)}\simeq
	\begin{cases}
		\{X_iX_j=0\}\subset\b P^3,&\kappa_{ij}<0,\\
		\{X_0X_{ij}=X_iX_j\}\subset\b P^3,&\kappa_{ij}=0,\\
		\{X_0X_{ij}=0\}\subset\b P^3,&\kappa_{ij}>0.
	\end{cases}
	\]
	Thus the fiber is $\b P^1\times\b P^1$ when $\kappa_{ij}=0$,
	and is the reduced union of two projective planes meeting along
	a line when $\kappa_{ij}\ne0$.
	
	For two points of the same color, there are homogeneous coordinates
	$[X_0:U:V:X_2]$ on $\b P^3$ in which
	\[
	Z^{\c H}_{2(p,i)}\simeq
	\begin{cases}
		\{V^2=0\}\subset\b P^3,&\kappa_{ii}<2,\\
		\{X_0X_2+V^2=0\}\subset\b P^3,&\kappa_{ii}=2,\\
		\{X_0X_2=0\}\subset\b P^3,&\kappa_{ii}>2.
	\end{cases}
	\]
	These are, respectively, a double plane, a quadric cone, and
	the reduced union of two projective planes meeting along a line.
\end{ex}

\subsubsection{Locality structure}\label{local}
Here we recall the definition of \textit{local spaces} in \textit{loc. cit.}, section 2.3.1. The locality structure on a space $Z$ over $\c H_{C\times I}$ is a system of isomorphisms $c_{D,E}$ for disjoint colored divisors $D,E$, compatible with commutativity and associativity:
$$c_{D,E}:\quad Z_{D}\times Z_E\xrightarrow{\sim} Z_{D\coprod E}.$$
We also ask $Z_\emptyset=\Spec\b C$. This is introduced in \cite{mirkovic2021loop}.

A local vector bundle $\c V\to \c H_{C\times I}$ is a vector bundle compatible with the locality structure: $\c V_D \otimes \c V_E\xrightarrow{\sim} \c V_{D+E}$.
\begin{rmk}
A Beilinson--Drinfeld \textit{factorization space} is a local space $Z\to \c H_{C\times I}$
such that the fiber $Z_D$ depends only on the formal neighborhood $\wh D$ of $D$ in $C$.\footnote{Equivalently $Z$ has a morphism to Ran space $\op{Ran}_C$.} For instance, Beilinson--Drinfeld Grassmannian is a factorization space over $\c H_{C\times I}$. We will reconstruct Beilinson--Drinfeld Grassmannian as a direct limit of Zastava spaces. 
\end{rmk}

\begin{rmk}
	In \cite{finkelberg1999note} and \cite{braverman2016gaiotto}, explicit coordinate systems are constructed for (open) Zastava spaces, which play an essential role in the study of their symplectic geometry and in applications to monopole moduli spaces. 
	Furthermore, the behavior of these coordinates under diagram automorphisms and folding is investigated in 	\cite{nakajima2023coulomb}.
	
	In contrast, the present work is concerned with \emph{compactified} Zastava spaces. The present construction does not use explicit (projective) coordinates.
\end{rmk}

\subsection{Interval Zastava space}
\subsubsection{Intervals}

We first generalize the notion of Hilbert scheme of divisors in section \ref{Hilbertdivisor} to arbitrary divisors.

For a pair $(D,E)\in \Div_{C\times I}^2$ with $D\le E$, i.e.\ $E-D$ is 
effective, we define the \emph{interval Hilbert scheme}
\[\c I(D,E) := D +\c H_{E-D}\subset\Div_{C\times I}.\]
This is a finite scheme containing all divisors $F$ satisfying $D\le F\le E$.

Let
$\Int := \{(D,E)\in \Div_{C\times I}^2 \mid D\le E\}$
be the space of pairs of divisors, here we understand it as a functor as $\Div_C$ is, i.e. a sheaf with respect to the flat topology. It is isomorphic to $\Div_{C\times I}\times \c H_{C\times I}$ by sending a pair $(D,E)$ to $(D,M=E-D)$.
Consider the following correspondence 
\[
\Int \xleftarrow{\,p\,} 
\c I := \{(D,E,F)\mid D\le F\le E\}
\xrightarrow{\,q\,}
\Div_{C\times I},
\]
where the projections are 
$p(D,E,F)=(D,E), q(D,E,F)=F$. $\c I$ is a closed subsheaf of $\Int\times \Div_{C\times I}$. (For the definition of closed subsheaf, one can see \cite[Definition 2.11]{stromme1996elementary} for instance.) For each $(D,E)\in \Int$, its fiber is $p^{-1}(D,E)=\c I(D,E)$, which is finite and $0$--dimensional.  
Similar to the case in \ref{c V}, $p$ is finite locally free of rank $2^{\deg(E-D)}$ \cite[Lemma 2.4.2]{mirkovic2021loop}.

Hence we can define $\b V := (p_* q^* \mathbb L)^*$ as a vector bundle over $\Int$, whose rank at $(D,E)$ is $2^{\deg(E-D)}$.
We will call $\b V$ the \emph{interval local vector bundle} coinduced by $\b L$, the locality for this interval setting will be explained in section \ref{subsec:interval-locality}.

For any pair $(D\le E)\in \Int$, using flat base change which is similar to the original case stated in the previous section, 
$\b V_{(D,E)}=\Gamma(\c I(D,E),\mathbb L)^*.$
In particular, $\b V_{(0,D)}=\Gamma(\c H_D,\b L^\kappa)^*=\c V_D$. Hence this  generalizes the bundle $\c V$ introduced in section \ref{recollection}.

\subsubsection{Zastava construction}
\label{intervalzastava}
Consider a pair of divisors $(D\le E)\in\Int$, and denote the difference by $M:=E-D\in\c H_{C\times I}.$
We define the \emph{left regular locus} $\Int^{\mathrm{reg,l}}\subset\Int$ consisting of all pairs $(D,E)$ such that
\begin{enumerate}
	\item $M=E-D$ is regular, i.e. its uncolored support consists of pairwise distinct points.
	\item $D$ and $M$ are disjoint.
\end{enumerate}
For a fixed colored divisor $A$, let
$\Int_A:=\{(A,B)\in\Int\mid A\leq B\}=A+\c H_{C\times I}$.  Then $\Int^{\mathrm{reg,l}}_A:=\Int_A\times_{\Int}\Int^{\mathrm{reg,l}}$ corresponds to the regular effective divisors $M\in\c H^{\op{reg}}$ such that its support
is disjoint from that of $A$. Thus $\Int_A^{\mathrm{reg,l}}$
is open and dense in $\Int_A$. We will use this observation in the proof of Theorem~\ref{thm:zastava-comparison}.

\begin{rmk}\label{leftrightdisjoint}
Similarly, define the right regular locus $\Int^{\rm{reg,r}}$ by replacing condition (2) with the condition that $E$ is disjoint from $M$. One can check that right regular locus $\Int^{\rm{reg,r}}$ does not intersect with the left one except $M=0$. However, we will show later in proposition \ref{prop:left-right-zastava} that the zastava spaces constructed in either way are canonically isomorphic, if $C$ has an open embedding to $\b A^1_z$.
\end{rmk}

For a pair $(D\le E)\in\Int^{\mathrm{reg,l}}$ with $M = E-D = \sum_{p\in S}p,$
the disjointness of $D$ and $M$ implies the locality isomorphism
$\c I(D,E)=D+\c H_M\cong D+\prod \c H_p$, and
\[
\b V_{(D,E)}=\Gamma(\c I(D,E),\mathbb L)^*
\cong
\mathbb L_D^* \otimes
\bigotimes_{p\in M} \b V_{(0,p)},
\quad \text{where }
\b V_{(0,p)}= \Gamma(\c H_p,\mathbb L)^*=\b L^*_0\oplus \b L^*_p.
\]

We now define the \emph{regular fibers} of the interval Zastava space.  
For $(D,E)\in\Int^{\mathrm{reg,l}}$ with $M=E-D=\sum_{p\in M}p$, set
\[
Z^{\mathrm{reg,l}}_{(D,E)}
:= \operatorname{im}\left(
\b P(\b L_D^*)\times \prod_{p\in M}\b P(\b V_{(0,p)})\hookrightarrow\b P(\b V_{(D,E)}) \right)\]
where the arrow is the Segre embedding induced by the above tensor product
decomposition of $\b V_{(D,E)}$. Note this construction depends only on the difference $E-D$ up to canonical isomorphism.

Finally, the \emph{interval Zastava space} is defined as the
scheme-theoretic closure
\[
Z:=\overline{Z^{\mathrm{reg,l}}}
\subset\mathbb P(\mathbb V).
\]
Theorem~\ref{thm:zastava-comparison} will show that this closure
restricts to the prescribed regular family over
$\Int^{\mathrm{reg,l}}$ and identify it with the pullback of
Mirkovi\'c's Zastava space under the difference map.

\subsubsection{Restriction to subsets of colors}
We now describe the behavior of interval Zastava spaces under
restriction of colors.
\begin{lem} \label{lem:based-product-decomposition}
	Suppose that the lattice $(\mathbb Z[I],\kappa)$ splits as a direct sum
	$(\mathbb Z[I_1],\kappa_1)\oplus(\mathbb Z[I_2],\kappa_2)$.
	Let $\Int,\Int_1$, and $\Int_2$ denote the corresponding spaces of intervals of $I, I_1, I_2$-colored divisors respectively.
	Then there is a canonical isomorphism
	$Z^\kappa\cong Z^{\kappa_1}\times Z^{\kappa_2}$.
	More precisely, let $D=D_1+D_2$ and $E=E_1+E_2$, where
	$D_i\leq E_i$ are $I_i$-colored divisors for $i=1,2$, then 
	\[
	Z^\kappa_{(D,E)}
	\cong
	Z^{\kappa_1}_{(D_1,E_1)}
	\times
	Z^{\kappa_2}_{(D_2,E_2)}.
	\]
	
	More generally, for any subset $I_1\subset I$, if $\kappa_1$ is the restriction of \(\kappa\) to \(\mathbb Z[I_1]\), then extension induces a canonical closed embedding \(Z^{\kappa_1}\hookrightarrow Z^\kappa\) over  \(\operatorname{Int}_1\hookrightarrow\operatorname{Int}\).
\end{lem}

\begin{proof}
	We only prove the product decomposition, the second assertion is similar.
	Let $\mathbb L_\kappa,\mathbb L_{\kappa_1}$, and
	$\mathbb L_{\kappa_2}$ be the corresponding Faltings line bundles.
	Since the two lattice summands are orthogonal with respect to $\kappa$, one has
	$\mathbb L_\kappa\simeq
	\mathbb L_{\kappa_1}\boxtimes\mathbb L_{\kappa_2}$.
	
	The decomposition of the colors also gives $\c I(D_1+D_2,E_1+E_2)
	\cong 	\c I(D_1,E_1)\times\c I(D_2,E_2).$
	Hence
	\begin{align*}
		\mathbb V^\kappa_{(D_1+D_2,E_1+E_2)}
		&=
		\Gamma\bigl(\c I(D_1+D_2,E_1+E_2),\mathbb L_\kappa\bigr)^*\\
		&\cong
		\Gamma\bigl(
		\c I(D_1,E_1)\times\c I(D_2,E_2),
		\mathbb L_{\kappa_1}\boxtimes\mathbb L_{\kappa_2}
		\bigr)^*\\
		&\cong
		\left(
		\Gamma\bigl(\c I(D_1,E_1),\mathbb L_{\kappa_1}\bigr)
		\otimes
		\Gamma\bigl(\c I(D_2,E_2),\mathbb L_{\kappa_2}\bigr)
		\right)^*\\
		&\cong
		\mathbb V^{\kappa_1}_{(D_1,E_1)}
		\otimes
		\mathbb V^{\kappa_2}_{(D_2,E_2)}.
	\end{align*}
	
	Under the identification $\Int\simeq\Int_1\times\Int_2$, the left
	regular locus $\Int^{\mathrm{reg,l}}$ consists of pairs of left regular
	intervals $(D_i\leq E_i)$ such that, after forgetting the colors,
	$E_1-D_1$ and $E_2-D_2$ have disjoint supports and
	$D_1+D_2$ is disjoint from
	$(E_1-D_1)+(E_2-D_2)$. 
	The Segre embedding $\operatorname{Seg}\colon
	\mathbb P(\mathbb V^{\kappa_1})
	\times
	\mathbb P(\mathbb V^{\kappa_2})
	\lhook\joinrel\longrightarrow
	\mathbb P(\mathbb V^\kappa)$ gives an isomorphism	
	\[
	Z^\kappa_{(D_1+D_2,E_1+E_2)}
	\cong
	Z^{\kappa_1}_{(D_1,E_1)}
	\times
	Z^{\kappa_2}_{(D_2,E_2)}.
	\]
\end{proof}
\begin{rmk}
	\label{rmk:affine-a-congruence}
	The product decomposition in the preceding lemma also holds
	for Mirkovi\'c Zastava spaces over colored Hilbert scheme: $Z^{\c H,\kappa}
	\simeq
	Z^{\c H,\kappa_1}\times Z^{\c H,\kappa_2}$. This provides a way to distinguish Zastava spaces associated with integrally congruent matrices.
	
	Let $\f g$ be a simple Lie algebra of type ADE, and let
	$\kappa$ and $\wt\kappa$ be the Cartan matrices of $\f g$
	and its untwisted affine Lie algebra respectively.
	Denote the finite simple roots by
	$\alpha_1,\ldots,\alpha_r$, the highest root by $\theta$,
	and the null imaginary root by $\delta$.
	The matrix $\wt\kappa$ is the Gram matrix in the basis
	$(\alpha_1,\ldots,\alpha_r,\delta-\theta)$,
	whereas $\kappa\oplus(0)$ is the Gram matrix in the basis
	$(\alpha_1,\ldots,\alpha_r,\delta).$
	So $\wt\kappa$ and $\kappa\oplus(0)$ are congruent.
	
	However, over the colored Hilbert scheme,
	\[
	Z^{\c H,\wt\kappa}
	\not\simeq
	Z^{\c H,\kappa\oplus(0)}
	\simeq
	Z^{\c H,\kappa}\times Z^{\c H,(0)}.
	\]	
	For $C=\b A^1$, these examples can be studied through
	Dong's identification of Mirkovi\'c Zastava spaces with
	compactified Coulomb branches \cite{dong2026loop}.
	The corresponding Coulomb branch descriptions can be found in
	\cite[\S3(vi), \S3(vii)]{BFN2}.
\end{rmk}

\subsection{Central extensions and the Faltings line bundle}

This subsection provides the technical preparation for
the translation isomorphisms of
Theorem~\ref{thm:interval-translation}.
We recall the central extensions of the lattice $L=\b Z[I]$
and the torus loop group $T_\c K$ associated with $\kappa$.
The latter extension defines a line bundle $\c R$ on $\Gr_T$.
We then compare its Beilinson--Drinfeld version with the
dual of the pullback of the Faltings line bundle.

\subsubsection{Central extension of lattice}\label{wh L}
Put $L=\b Z[I]$ and
$T=L\otimes_{\b Z}\b G_m$, and recall that $\kappa$ is even and integral.
Fix once and for all a 2-cocycle
$\epsilon:L\times L\to\b Z/2\b Z$ satisfying
\[
\epsilon(\lambda,0)=\epsilon(0,\lambda)=0,\qquad
\epsilon(\lambda,\mu)+\epsilon(\mu,\lambda)
\equiv\kappa(\lambda,\mu)\pmod 2.
\]
Such a cocycle exists, for instance with basis $\{e_i\}$ of $L$ one can take
\[
\epsilon\left(\sum_i n_i e_i,\sum_j m_j e_j\right)
=\sum_{i>j}n_i m_j\kappa(e_i,e_j)\pmod 2.
\]
The 2-cocycle $\epsilon$ defines a central extension
\[
1\longrightarrow\b G_m\longrightarrow\wh L
\longrightarrow L\longrightarrow0,
\]
with multiplication
$(a,\lambda)(b,\mu)
=\bigl(ab(-1)^{\epsilon(\lambda,\mu)},\lambda+\mu\bigr)$, $a,b\in \b G_m, \lambda,\mu\in L$.
Thus lifts of $\lambda,\mu$ have commutator $(-1)^{\kappa(\lambda,\mu)}$.

\subsubsection{Central extension of $T_\c K$}
\label{subsubsec:torus-central-extension}
We recall the construction of \cite[3.2.1]{zhu2009affine},
\cite[17.1]{GaitsgoryCFT}, and
\cite[1.4--1.7]{BeilinsonHeisenberg}. 
The triple $(\b Z[I],\kappa,\epsilon)$ determines a symmetric (Heisenberg) central extension
\[
1\longrightarrow\b G_m\longrightarrow\wh T_\kappa
\xrightarrow{p}T_{\c K}\longrightarrow1
\]
with a canonical splitting
$i:T_{\c O}\to\wh T_\kappa$. Its commutator pairing is $[\lambda(f),\mu(g)]
=\{f,g\}^{-\kappa(\lambda,\mu)},$
where $\{f,g\}\in\b G_m$ is the Contou--Carr\`ere symbol, see \cite[3.2.1]{zhu2009affine}. Moreover, there is a section $s:T_\c K\to \wh T_\kappa$, such that for any $\lambda,\mu\in L$, $s(z^\lambda)s(z^\mu)=(-1)^{\epsilon(\lambda,\mu)}s(z^{\lambda+\mu})$. (For instance one can use the arguments in \cite[5.5]{kac1998vertex}.)

The homomorphism $L\to T_{\c K}$ sending $\lambda$ to $z^\lambda$ defines a fiber product 
$\wh T_\kappa\times_{T_{\c K}}L$.
One can check there is a group isomorphism from the central extension $\wh L$ introduced in \ref{wh L} to this fiber product 
$\wh T_\kappa\times_{T_{\c K}}L$, sending 
$(1,\lambda)$ to $(s(z^\lambda),\lambda)$.

The canonical splitting $i$ defines a $\b G_m$-torsor
\[
\wh T_\kappa/i(T_{\c O})
\longrightarrow T_{\c K}/T_{\c O}=\Gr_T.
\]
We denote the associated line bundle on $\Gr_T$ by $\c R$.\footnote{Our
$\c R$ is written as $R_Q$ in \cite[17.1]{GaitsgoryCFT} and
as $\c L_T^{-1}$ in \cite[3.2.1]{zhu2009affine}.}

From the isomorphism $\wh L\simeq \wh T_\kappa\times_{T_{\c K}}L$, for any $\lambda\in L$ we can identify the nonzero vectors in the
corresponding fiber:
$\wh L_\lambda
\xrightarrow{\sim}\c R^\times_{z^\lambda}$, sending $
(a,\lambda)$ to $[a\cdot s(z^\lambda)]:=as(z^\lambda)T_\c O.$

 Over the curve $C$, there is also a Beilinson--Drinfeld version of central extension, see \cite[2.2--2.4]{BeilinsonHeisenberg}:
\[
1\longrightarrow\b G_{m,C}\longrightarrow\wh T_{\kappa,C}
\longrightarrow T_{\c K,C}\longrightarrow1.
\]
with a canonical splitting
$i_C:T_{\c O,C}\to\wh T_{\kappa,C}.$
The quotient
$\c R_C^\times=\wh T_{\kappa,C}/i_C(T_{\c O,C})$
defines a line bundle $\c R_C$ on
$\Gr_{T,C}=T_{\c K,C}/T_{\c O,C}$. Similarly by considering the central extension over $C^n$ we can define a line bundle
$\c R_{C^n}$ on $\Gr_{T,C^n}$ compatible with factorization.

\subsubsection{Comparison of line bundles $\c R$ and Faltings line bundle $\b L$}
\label{Comparison of R and L}

Recall an $R$-point of the Beilinson--Drinfeld
Grassmannian $\Gr_{T,C^n}$ can be described as the moduli of tuples 
$(x_1,\ldots,x_n,P,s)$, where $x_a:S=\Spec R\to C$, $P$ is a $T$-bundle over $C$, and $s\in\Gamma(C-\cup_i x_i,P)$ is a section. 

Take $D=(D_i)_{i\in I}$ to be the negative of corresponding colored divisor of section $s$, then $D$ gives the corresponding $T$-bundle 
$(\c O_{C\times S}(-D_i))_{i\in I}$ together with canonical section (equivalently, trivialization) away from these graphs. As a result,
 $\Gr_{T,C^n}$ can be also described as the moduli of tuples 
$(x_1,\ldots,x_n,D)$, where $x_a:S\to C$ and
$D=(D_i)_{i\in I}$ is a relative colored Cartier divisor supported on the union of the graphs of
the $x_a$. We will use this alternative description in this subsection and the next.

Forgetting the marked points defines a morphism
\[
\rho_n:\Gr_{T,C^n}\longrightarrow\Div_{C\times I},
\qquad
(x_1,\ldots,x_n,D)\longmapsto D.
\]

\begin{prop}
	\label{prop:faltings-heisenberg-comparison}
	For every $n$, there is a canonical isomorphism 
	\[
	\rho_n^*\b L\xrightarrow{\sim}\c R_{C^n}^{-1},
	\]
	compatible with factorization structures.
\end{prop}
\begin{proof}
	Write \(j_\gamma:C\to \Gr_{T,C}\) sending $x$ to the element corresponding to \(\gamma[x]\), and \(j_{\gamma,\delta}: C^2\to \Gr_{T,C^2}\) sending $(x,y)$ to the element corresponding to \(\gamma[x]+\delta[y]\). This is slightly different from Remark~\ref{thetadata}.	By \cite[Lemma~17.1]{GaitsgoryCFT},
	\[
	\lambda_{\c R}^{\gamma}:=j_\gamma^*\c R_C
	\simeq \omega_C^{-\kappa(\gamma,\gamma)/2}.
	\]
	For two points collision, \cite[Lemma~17.3]{GaitsgoryCFT} gives
	\[
	j_{\gamma,\delta}^*\c R_{C^2}
	\simeq
	(\lambda_{\c R}^{\gamma}\boxtimes\lambda_{\c R}^{\delta})
	\otimes\c O_{C^2}\bigl(\kappa(\gamma,\delta)\Delta\bigr).
	\]
	Restricted on diagonal, $j_{\gamma,\delta}^*\c R_{C^2}|_{\Delta}$ gives line bundle	\(\lambda_{\c R}^{\gamma+\delta}\).
	 Since
	\(\c O_{C^2}(-\Delta)|_\Delta\simeq\omega_C\), restriction of this
	isomorphism gives the collision map
	\[
	c_{\c R}^{\gamma,\delta}\colon
	\lambda_{\c R}^{\gamma}\otimes\lambda_{\c R}^{\delta}
	\xrightarrow{\sim}
	\lambda_{\c R}^{\gamma+\delta}\otimes\omega_C^{\kappa(\gamma,\delta)}.
	\]
	The $\theta$-data of $\b L$ is explained in	Remark~\ref{thetadata},  its pullback $\rho_n^*\b L$ has the same $\theta$-data with $\c R_{C^n}^{-1}$, as factorizable line bundles over $\Gr_T$. By \cite[Lemma~1.3 and Proposition~1.4]{tao2021central}, we get the required isomorphisms for all \(n\).
\end{proof}
\begin{rmk}
	After completing this work, we notice
	\cite[Theorems~6.1--6.3]{osipov2026contou}, where Osipov
	explains the relation between Deligne pairings and the
	Contou--Carr\`ere symbol.
\end{rmk}

\subsection{Comparison of interval Zastava spaces}
We first construct translation isomorphisms for the
projective bundle $\b P(\b V)$.
We then use these isomorphisms to compare interval
Zastava spaces with Mirkovi\'c's construction.
We establish the locality structure for interval Zastava spaces.
We also compare the versions of interval Zastava spaces
defined using the left and right regular loci.

\subsubsection{Translation isomorphisms}
We now compare the fibers of the projectivization of the factorization vector bundle $\b V$. For the rest of this subsection, we assume that $C$ is an open subset of $\b A^1_z$, where $z$ is a global \'etale coordinate.\footnote{When defining Deligne pairing and Faltings line bundle, we ask the curve to be smooth projective.
So here we are really treating $C$ as an open subset of the curve $\b P^1\supset \b A^1_z$.\label{footnote}}

\begin{thm}
	\label{thm:interval-translation}
	The chosen coordinate $z$ determines canonical
	isomorphisms
	\[
	f_D:\b P\b V_{(A,B)}
	\xrightarrow{\sim}
	\b P\b V_{(A+D,B+D)}
	\]
	for colored divisors $A\leq B$ and $D$.
	 They satisfy
	$f_E\circ f_D=f_{D+E}$.
\end{thm}

The strategy is to use the description
$\b V_{(A,B)}=\Gamma(\c I(A,B),\b L)^\vee$
and compare these spaces for translated intervals.
Proposition~\ref{prop:faltings-heisenberg-comparison}
identifies the pullback of $\b L$ to the
line bundle $\c R^{-1}_{C^n}$.
Translation (i.e. addition) of divisors is realized by an element in the central extension $\wh T_\kappa$ acting on $\c R$.
This action induces a comparison of the spaces of
sections of $\b L$, and we can descend it to the projectivizations.

\begin{proof}
	
	We work locally over a test scheme $S=\Spec R$.
 Addition of $D$ defines an isomorphism
	\[
	\tau_D:\c I(A,B)
	\xrightarrow{\sim}
	\c I(A+D,B+D),
	\qquad F\longmapsto F+D.
	\]
	
	Choose a map $\mathbf x=(x_1,\ldots,x_n):S\to C^n$
	whose graphs contain the supports of the divisors $A,B,D$. Consider morphisms
	\[
	\begin{aligned}
		\eta:\c I(A,B)&\longrightarrow\Gr_{T,C^n},
		&F&\longmapsto(\mathbf x,F),\\
		\eta_D:\c I(A+D,B+D)
		&\longrightarrow\Gr_{T,C^n},
		&F'&\longmapsto(\mathbf x,F').
	\end{aligned}
	\]
		Proposition~\ref{prop:faltings-heisenberg-comparison}
	gives
$$\eta^*\c R_{C^n}^{-1}
\simeq \eta^*\rho_n^*\b L
\simeq \b L|_{\c I(A,B)},\quad 
\eta_D	^*\c R_{C^n}^{-1}
\simeq \eta_D^*\rho_n^*\b L
\simeq \b L|_{\c I(A+D,B+D)}.$$
	
	For each color $i$, write $D_i=D_i^+-D_i^-$ and let
	$f_i^\pm$ be the canonical monic defining polynomials of $D_i^\pm$.
	The ratio defines an element in $T_{\c K,C^n}$,  the Beilinson--Drinfeld
	torus loop group over $C^n$: 
	\[
	t_D
	:=\left(\frac{f_i^+}{f_i^-}\right)_{i\in I}.
	\]
	Let $t_D$ also denote its left action on $\Gr_{T,C^n}$. Hence $\eta_{D}\circ\tau_D
	=t_D\cdot\eta.$
	
Choose a lift $\widetilde t_D$ of $t_D$ in $\wh T_{\kappa,C^n}$, and let
$\widetilde t_D$ itself also denote its action on $\c R_{C^n}$
induced by left multiplication sending $[\widehat g]$
to $[\widetilde t_D\widehat g]$.
Hence multiplication by $\widetilde t_D$ defines a
map of line bundles
\[
\widetilde\tau_D:
\eta^*\c R_{C^n}
\xrightarrow{\sim}
\eta_D^*\c R_{C^n},
\qquad
(F,v\in \left( \eta^*\c R_{C^n})_F\right) \mapsto
(F+D,\widetilde t_D\cdot v).
\]

In terms of line bundle $\b L$, $\widetilde\tau_D$ gives $\b L^*|_{\c I(A,B)}\to \b L^*|_{\c I(A+D,B+D)}$. Dualizing and then take global sections we get
$$\Gamma\bigl(\c I(A+D,B+D),\b L\bigr)
\xrightarrow{\sim}
\Gamma\bigl(\c I(A,B),\b L\bigr).$$
Dualizing again, $\wt \tau_D$ induces a map
\[
\Phi_{\widetilde t_D}:
\b V_{(A,B)}
=\Gamma\bigl(\c I(A,B),\b L\bigr)^\vee
\xrightarrow{\sim}
\Gamma\bigl(\c I(A+D,B+D),\b L\bigr)^\vee
=\b V_{(A+D,B+D)}.
\]

Replacing $\widetilde t_D$ by $c\widetilde t_D$,
	where $c\in R^\times$, multiplies
	$\Phi_{\widetilde t_D}$ by $c$.
	Its projectivization is therefore independent of the
	lift. Hence $\Phi_{\widetilde t_D}$ descends to a canonical isomorphism
	\[
	f_D:\b P\b V_{(A,B)}
	\xrightarrow{\sim}
	\b P\b V_{(A+D,B+D)}.
	\]	
	
Finally, for any $D,E$, the fact $t_Et_D=t_{D+E}$
and independence of the lift after projectivization imply
\[	f_E\circ f_D=f_{D+E}.\]
\end{proof}

\begin{rmk}
	The construction uses the chosen open embedding
	$C\hookrightarrow\b A^1_z$ to obtain canonical rational functions representing the divisors.
	An arbitrary \'etale morphism $C\to\b A^1$ does not supply such representatives, and the theorem is not true in this general case.
\end{rmk}

\subsubsection{Comparison with Mirkovi\'c's construction}
Continue with the chosen open embedding $C\hookrightarrow\b A^1_z$.
Denote the difference map by
\[
\delta:\Int\longrightarrow\c H_{C\times I},\qquad (D,E)\longmapsto E-D.
\]
Let $Z^{\c H}\subset\b P(\c V)$ be the Mirkovi\'c Zastava space over $\c H_{C\times I}$ recalled in section \ref{recollection}, and let $Z\subset\b P(\b V)$ be the interval Zastava space over $\Int$. 
By Theorem~\ref{thm:interval-translation}, translation gives an
isomorphism of projective bundles over $\Int$,
\[
\begin{gathered}
	f:\delta^*\b P(\c V)\xrightarrow{\sim}\b P(\b V),\\
	f|_{(D,E)}=f_D:
	\b P(\c V_{E-D})\xrightarrow{\sim}\b P(\b V_{(D,E)}).
\end{gathered}
\]
Here we have $\c V_M=\b V_{(0,M)}$.

\begin{thm}
	\label{thm:zastava-comparison}
	The isomorphism $f$ restricts to an isomorphism over $\Int$,
	\[
	f:\delta^*Z^{\c H}\xrightarrow{\sim}Z,
	\]

As a corollary, for every interval $A\leq B$ and every colored divisor $D$, translation
restricts to an isomorphism
\[
f_D:Z_{(A,B)}\xrightarrow{\sim}Z_{(A+D,B+D)}.
\]
\end{thm}

\begin{proof}
	By definition of $Z^{\mathrm{reg,l}}$, $	\left.f(\delta^*Z^{\c H})\right|_{\Int^{\mathrm{reg,l}}}
		=Z^{\mathrm{reg,l}}.$
		Since $f(\delta^*Z^{\c H})$ is closed, the definition of $Z$
		as a scheme-theoretic closure gives
		$Z\subseteq f(\delta^*Z^{\c H})$.
		
		To prove the reverse inclusion, first fix a colored divisor $A$.
		Under the identification
		$\phi:\c H_{C\times I}\xrightarrow{\sim}\Int_A$, the restriction
		$f^{-1}(Z)|_{\Int_A}$ is a closed subscheme of $Z^{\c H}$
		containing the restriction of $Z^{\c H}$ to $\phi^{-1}(\Int_A^{\mathrm{reg,l}})$.
		Since $\Int_A^{\mathrm{reg,l}}$ is open and dense in $\Int_A$, and 
		the restriction of $Z^{\c H}$ to $\phi^{-1}(\Int_A^{\mathrm{reg,l}})$ is dense in $Z^{\c H}$.
		Moreover, $Z^{\c H}$ is reduced, as it is the scheme-theoretic
		closure of its reduced regular family.
		A closed subscheme of a reduced scheme containing a dense
		open subscheme is the whole scheme. Hence
		\[
		\left.f^{-1}(Z)\right|_{\Int_A}=Z^{\c H}.
		\]
		
	To obtain the equality in families, we need to locally take a presentation of interval in terms of effective divisors. Consider the map
	\[
	\psi:(\c H_{C\times I})^3\longrightarrow\Int,\qquad
	(P,Q,M)\longmapsto(P-Q,P-Q+M).
	\]
	The pullback $\psi^*f^{-1}(Z)$ is a closed subscheme of
	$\psi^*\delta^*Z^{\c H}	\simeq(\c H_{C\times I})^2\times Z^{\c H}$.
	For every closed point $(P,Q)$, the comparison on $\Int_{P-Q}$ identifies its fiber with the whole corresponding $Z^{\c H}$. Since $(\c H_{C\times I})^2\times Z^{\c H}$ is reduced, $\psi^*f^{-1}(Z)$ equals $(\c H_{C\times I})^2\times Z^{\c H}$. Since the local presentation $\psi$ is arbitrary, we get $f^{-1}(Z)=\delta^*Z^{\c H}$.
	\end{proof}

\subsubsection{Locality structure}
\label{subsec:interval-locality}

We extend the notion of locality from Section~\ref{local}
to interval spaces.

For an interval $(D,E)$, set
$\supp(D,E):=\bigcup_{i\in I}\bigl(\supp D_i\cup\supp E_i\bigr).$
Two intervals are called \emph{disjoint} if their supports
are disjoint.
A \emph{locality structure} on a space $Y\to\Int$ consists
of isomorphisms
\[
c_{(D_1,E_1),(D_2,E_2)}:
Y_{(D_1,E_1)}\times Y_{(D_2,E_2)}
\xrightarrow{\sim}
Y_{(D_1+D_2,E_1+E_2)}
\]
for disjoint intervals over $S$, compatible with base
change, associativity, commutativity, and $Y_{(0,0)}\simeq\Spec\b C$.
For vector bundles, the analogous definition uses tensor
products.

\begin{lem}
	\label{prop:interval-vector-locality}
	The factorization structure of $\b L$ in
	Proposition~\ref{locality} induces a canonical locality
	structure on the vector bundle $\b V$ over $\Int$.
	\qed
\end{lem}

\begin{cor}
	\label{cor:interval-locality}
	The comparison isomorphism	of Theorem~\ref{thm:zastava-comparison} transports the 	locality structure of Mirkovi\'c Zastava to $Z$.\qed	
\end{cor}

\begin{rmk}
	For a general curve, locality of $Z$ does not
	follow directly from its definition by closure.
	The left regular locus $\Int^{\mathrm{reg,l}}$ is dense
	but not open in $\Int$, so taking closure need not
	preserve the prescribed regular fibers.
\end{rmk}

\subsubsection{Comparison of the left and right interval constructions}
 Continue to fix the embedding $C\hookrightarrow\b A^1_z$ . Let $Z^L$ and $Z^R$ be the
scheme-theoretic closures in $\b P(\b V)$ of the regular families over
$\Int^{\reg,l}$ and $\Int^{\reg,r}$, respectively.

\begin{prop}
	\label{prop:left-right-zastava}
	The choice of embedding $C\hookrightarrow\b A^1_z$ gives a canonical isomorphism $Z^L\simeq Z^R$ over $\Int$.
\end{prop}

\begin{proof}
	For every interval $(D,E)$, negating divisors gives an isomorphism
	$\c I(D,E)\xrightarrow{\sim}\c I(-E,-D)$ by sending $F$ to $-F$.
	By Lemma~\ref{dual}, the canonical isomorphisms
	$\b L_F\simeq\b L_{-F}$ induce an isomorphism
	\[
	\iota_{D,E}:\b P\b V_{(D,E)}
	\xrightarrow{\sim}\b P\b V_{(-E,-D)}.
	\]
	These isomorphisms assemble into an isomorphism of the ambient projective
	bundles covering the involution $(D,E)\mapsto(-E,-D)$ of $\Int$.
	
	Recall $(D,E)\in\Int^{\op{reg,r}}$ means $E$ and $M=E-D$ are
	disjoint, and $M$ is regular in $\c H_{C\times I}$. 
	Thus $(-E,-D)$ lies in $\Int^{\reg,l}$. $\iota$ induces an isomorphism on $(Z^L)^{\text{reg,l}}$ and $(Z^R)^{\text{reg,r}}$, hence taking closures gives
	\[	\iota_{D,E}\bigl(Z^R_{(D,E)}\bigr)	=Z^L_{(-E,-D)}.	\]
	
	Applying the translation isomorphism in Theorem~\ref{thm:zastava-comparison}, $Z^L_{(-E,-D)}$ is isomorphic to $Z^L_{(D,E)}$.
\end{proof}

\subsection{Interval Zastava spaces and affine Grassmannians}
\label{subsec:interval-grassmannian}

Throughout this subsection, keep the chosen open embedding
$C\hookrightarrow\b A^1_z\subset\b P^1$.
We first compare punctual interval Zastava spaces with
intersections of semi-infinite orbit closures in the
affine Grassmannian.
We then construct an inductive system of interval
Zastava spaces and extend it to a family over $\c H_C$.
The construction applies to any based lattice
$(\b Z[I],\kappa)$. In the group case, it recovers the
affine Grassmannian and its Beilinson--Drinfeld version. 

\subsubsection{The group case}
\label{zastavagroup}

Assume that $G$ is a simple simply connected group of
type ADE.
Identify $\b Z[I]$ with $X_*(T)$ by sending the basis
vectors to the simple coweights, and let $\kappa$ be the
Cartan matrix. Choose opposite Borel subgroups $B=TN$ and $B^-=TN^-$.

Fix a point $p\in C$ and put $t=z-z(p)$.
For $\nu\in X_*(T)$, let
\[
S_{\nu,p}=N(\c K_p)t^\nu G(\c O_p)/G(\c O_p),
\qquad
T_{\nu,p}=N^-(\c K_p)t^\nu G(\c O_p)/G(\c O_p)
\]
be the corresponding semi-infinite orbits of $\Gr_{G,p}:=G(\c K_p)/G(\c O_p)\simeq \Gr_G$.

The orbit closures are
\[
\ol S_{\nu,p}=\bigcup_{\eta\leq\nu}S_{\eta,p},
\qquad
\ol T_{\nu,p}=\bigcup_{\eta\geq\nu}T_{\eta,p}.
\]
Moreover, $\ol S_{\nu,p}\cap\ol T_{\eta,p}$ is empty
unless $\eta\leq\nu$, i.e. $\nu-\eta$ is a
nonnegative integral combination of simple coweights. See \cite{mirkovic2007geometric}.
Multiplication by $t^\xi$ sends $S_{\nu,p}$ to
$S_{\nu+\xi,p}$ and $T_{\nu,p}$ to $T_{\nu+\xi,p}$.

For $\alpha\in\b N[I]$, the comparison theorem
of \cite[\S2.4.6]{mirkovic2021loop} identifies Mirkovi\'c zastava with intersections of closures of semiinfinite orbits:
\[
Z^{\c H}_{\alpha[p]}
\simeq
\ol S_{0,p}\cap\ol T_{-\alpha,p}.
\]
Combining this identification with
Theorem~\ref{thm:zastava-comparison} gives the
following interval version.

\begin{prop}
	\label{prop:interval-zastava-grassmannian}
	For $\lambda\leq\mu$, there is a canonical isomorphism
	\[
	\phi_{\lambda,\mu}:
	Z_{(\lambda[p],\mu[p])}
	\xrightarrow{\sim}
	\ol S_{-\lambda,p}\cap\ol T_{-\mu,p}.
	\]
\end{prop}

\begin{proof}
	Put $\alpha=\mu-\lambda$.
	By Theorem~\ref{thm:zastava-comparison},
	\[
	Z_{(0,\alpha[p])}
	\simeq
	Z^{\c H}_{\alpha[p]}
	\simeq
	\ol S_{0,p}\cap\ol T_{-\alpha,p}.
	\]
	Recall that $f_{\lambda[p]}$ is the translation
	isomorphism of Theorem~\ref{thm:interval-translation},
	whose restriction to Zastava spaces is established
	in Theorem~\ref{thm:zastava-comparison}.
	Define $\phi_{\lambda,\mu}$ by the commutative diagram
	\[
	\begin{tikzcd}
		Z_{(0,\alpha[p])}
		\arrow[r,"\sim"]
		\arrow[d,"f_{\lambda[p]}"']
		&
		\ol S_{0,p}\cap\ol T_{-\alpha,p}
		\arrow[d,"t^{-\lambda}"]
		\\
		Z_{(\lambda[p],\mu[p])}
		\arrow[r,"\phi_{\lambda,\mu}"]
		&
		\ol S_{-\lambda,p}\cap\ol T_{-\mu,p}.
	\end{tikzcd}
	\]
	Both vertical maps are isomorphisms.
\end{proof}

\subsubsection{The inductive system of interval Zastava spaces}
\label{subsec:interval-inductive-system}

We now return to an arbitrary even integral symmetric
form $\kappa$ on $\b Z[I]$.

Suppose that $F\leq D\leq E\leq H.$ The inclusion of interval Hilbert schemes $\c I(D,E)\hookrightarrow\c I(F,H)$
induces a surjection $\Gamma(\c I(F,H),\b L)
\twoheadrightarrow
\Gamma(\c I(D,E),\b L),$
and hence an inclusion
$\b V_{(D,E)}\hookrightarrow\b V_{(F,H)}.$

\begin{lem}
	\label{lem:interval-transition}
	Fix an open embedding $C\hookrightarrow\b A^1_z$.
	For colored divisors $F\leq D\leq E\leq H$, the
	canonical linear inclusion
	$\b V_{(D,E)}\hookrightarrow\b V_{(F,H)}$
	induces a closed immersion
	\[
	i_{(D,E),(F,H)}:
	Z_{(D,E)}\hookrightarrow Z_{(F,H)}.
	\]
	These maps are compatible with compositions.
\end{lem}

\begin{proof}
	It suffices to prove separately the cases $D=F$ and
	$E=G$, since a general inclusion factors as
	\[
	(D,E)\subset(F,E)\subset(F,G).
	\]
	
	First suppose that $D=F$.
	The common translation by $-D$ identifies the
	ambient projective embedding with
	\[
	\b P(\c V_{E-D})\hookrightarrow\b P(\c V_{G-D}).
	\]
	Indeed, the line bundle maps defining translation
	commute with restriction of sections to smaller
	interval Hilbert schemes.
	By Theorem~\ref{thm:zastava-comparison}, the
	corresponding Zastava spaces are identified with
	$Z^{\c H}_{E-D}$ and $Z^{\c H}_{G-D}$.
	The desired closed immersion is therefore the
	growth map
	\[
	Z^{\c H}_{E-D}\hookrightarrow Z^{\c H}_{G-D}
	\]
	of \cite[\S2.3.1, Lemma~(a)]{mirkovic2023lattice}.
	Since this map is also constructed by restriction
	of sections, it agrees with the given ambient
	projective embedding.
	
	To treat the second case, use divisor inversion
	$K\mapsto-K$.
	By Lemma~\ref{dual}, it induces isomorphisms
	\[
	\iota_{D,E}:
	\b P(\b V_{(D,E)})
	\xrightarrow{\sim}
	\b P(\b V_{(-E,-D)}).
	\]
	Their compatibility with restriction of sections
	gives the commutative diagram
	\[
	\begin{tikzcd}
		\b P(\b V_{(D,E)})
		\arrow[r,hook]
		\arrow[d,"\iota_{D,E}"']
		&
		\b P(\b V_{(F,G)})
		\arrow[d,"\iota_{F,G}"]
		\\
		\b P(\b V_{(-E,-D)})
		\arrow[r,hook]
		&
		\b P(\b V_{(-G,-F)}).
	\end{tikzcd}
	\]
	
	These inversion maps restrict to isomorphisms
	\[
	\iota_{D,E}:Z_{(D,E)}
	\xrightarrow{\sim}Z_{(-E,-D)}.
	\]
	Indeed, putting $M=E-D$, the comparison
	isomorphisms identify inversion with endpoint
	reversal $K\mapsto M-K$ on $\c H_M$.
	Over the regular effective locus, factorization
	identifies this reversal with automorphisms of
	the individual projective-line factors.
	It therefore preserves the regular Segre family
	and hence its closure $Z^{\c H}$.
	
	Now suppose that $E=G$.
	The inverted intervals
	\[
	(-E,-D)\subset(-E,-F)
	\]
	have the same lower endpoint, so the first case
	gives a closed immersion
	\[
	Z_{(-E,-D)}\hookrightarrow Z_{(-E,-F)}.
	\]
	Transporting this map by inversion gives
	$Z_{(D,E)}\hookrightarrow Z_{(F,E)}$.
	The commutative diagram above shows that this
	is precisely the restriction of the given
	ambient projective embedding.
\end{proof}

These maps are compatible with compositions and make the
interval Zastava spaces into a directed system, ordered by inclusion of intervals.

\subsubsection{Generalization of affine Grassmannian}
Let $\kappa$ be an integral symmetric matrix with even
diagonal entries.

For an effective divisor $A\in\c H_C$, define
\[\c G^\kappa_A
:=
\varinjlim_{\substack{	D\leq E\\\supp(D),\supp(E)\subset A}}
Z_{(D,E)}\subset\c G^\kappa,\]
here for colored divisor $D=(D_i)_{i\in I}$, $\supp(D)$ means the union of $\supp(D_i)$.
We know from definition $\c G^\kappa_A$ is ind-projective, and
replacing $A$ by a positive multiple does not change it.

Define a family $\c G^\kappa_{\c H_C}\to\c H_C$ whose fiber along $A\in \c H_C$ is $\c G^\kappa_A$. When $A$ and $B$ are disjoint effective divisors, Corollary \ref{cor:interval-locality} gives canonical isomorphisms
\[
\c G^\kappa_A\times\c G^\kappa_B
\xrightarrow{\sim}
\c G^\kappa_{A\coprod B}.
\]
These isomorphisms give
$\c G^\kappa_{\c H_C}$ the structure of a factorization space.

Define \textit{semiinfinite spaces} associated to $\kappa$, for any $D,E\in\Div_{C\times I}$:
\[
S_{-D}^\kappa
:=
\varinjlim_{\substack{E'\geq D\\\supp(E')\subseteq\supp(D)}}
Z_{(D,E')},
\qquad
T_{-E}^\kappa
:=
\varinjlim_{\substack{D'\leq E\\\supp(D')\subseteq\supp(E)}}
Z_{(D',E)}.
\]
They are also indschemes, having structures of factorization spaces.
Then we have
\[
\c G^\kappa_A
\simeq
\varinjlim_{D,\supp(D)\subset A}S^\kappa_{-D}
\simeq
\varinjlim_{E,\supp(E)\subset A}T^\kappa_{-E}.
\]

When $A=p\in C$ is a single point, 
\[
S^\kappa_{-\lambda[p]}
=
\varinjlim_{\mu\geq\lambda}
Z_{(\lambda[p],\mu[p])},
\qquad
T^\kappa_{-\mu[p]}
=
\varinjlim_{\lambda\leq\mu}
Z_{(\lambda[p],\mu[p])}.
\]

In the group case of section~\ref{zastavagroup},
Proposition~\ref{prop:interval-zastava-grassmannian}
identifies $S^\kappa_{-\lambda,p}, T^\kappa_{-\mu,p}$ with intersections of semi-infinite orbit closures:
\[
S^\kappa_{-\lambda[p]}\simeq\ol S_{-\lambda,p},
\qquad
T^\kappa_{-\mu[p]}\simeq\ol T_{-\mu,p}.
\qquad
\c G^\kappa_p\simeq\Gr_{G,p},
\qquad
\c G^\kappa_{\c H_C}
\simeq
\Gr^{\mathrm{BD}}_{G,\c H_C}.
\]
Thus $\c G^\kappa_p, \c G^\kappa_{\c H_C}$ generalize the notion of affine Grassmannian and Beilinson--Drinfeld Grassmannian 
to integral symmetric matrices with even diagonal entries.

\subsubsection{The line bundle $\c O(1)$ over $\c G^\kappa$} 
\label{subsec:grassmannian-line-bundle}

For every interval $D\leq E$, set
\[
\c O_{Z_{(D,E)}}(1)
:=
\left.
\c O_{\b P(\b V_{(D,E)})}(1)
\right|_{Z_{(D,E)}}.
\]
For $F\leq D\leq E\leq H$, the closed immersion
$i_{(D,E),(F,H)}: Z_{(D,E)}\hookrightarrow Z_{(F,H)}$
of Lemma~\ref{lem:interval-transition} is induced by
the  inclusion $\b V_{(D,E)}\hookrightarrow\b V_{(F,H)}$.
Hence
\[
i_{(D,E),(F,H)}^*
\c O_{Z_{(F,H)}}(1)
\simeq
\c O_{Z_{(D,E)}}(1).
\]
These isomorphisms are compatible with compositions.

Thus using the line bundles $\c O_{Z_{(D,E)}}(1)$, we can define a line
bundle $\c O(1)$ on the indscheme $\c G^\kappa_A$ and on
$\c G^\kappa_{\c H_C}$. It is compatible with locality:
for disjoint effective divisors $A$ and $B$, the
locality isomorphism
\[
c_{A,B}:
\c G^\kappa_A\times\c G^\kappa_B
\xrightarrow{\sim}
\c G^\kappa_{A\coprod B}
\]
gives
\[
c_{A,B}^*\c O_{\c G^\kappa_{A\coprod B}}(1)
\simeq
\c O_{\c G^\kappa_A}(1)
\boxtimes
\c O_{\c G^\kappa_B}(1).
\]

\subsubsection{Global sections and the lattice vertex algebra}
\label{subsec:lattice-vertex-algebra}

Let $\kappa$ be an integral symmetric matrix with even
diagonal entries, and let $L_\kappa$ denote the
associated lattice vertex algebra, see
\cite{frenkel2004vertex,kac1998vertex} for references.
We identify the dual global sections of
$\c O_{\c G^\kappa_p}(1)$ with $L_\kappa$, whose geometric form uses
$\Gr_{T,p}$, see \cite[\S3]{zhu2009affine}.
The proof uses the Zastava comparison theorem established earlier.

Keep the chosen open embedding
$C\hookrightarrow\b A^1_z$, and fix $p\in C$.
For an effective colored divisor $M$, by the construction of $Z^\c H$, the Kodaira embedding
\[
\c H_M\hookrightarrow\b P\Gamma(\c H_M,\b L)^*=\b P(\c V_M)
\]
factors through $Z^{\c H}_M$. See \cite[\S2.4.3]{mirkovic2021loop}.

For an interval $D\leq E$, put $M=E-D$. By the construction in
Theorem~\ref{thm:interval-translation}, the projective isomorphism
\[
f_D:
\b P(\c V_M)\xrightarrow{\sim}
\b P(\b V_{(D,E)})
\]
carries the Kodaira image of $\c H_M$ to that of $\c I(D,E)$.
Theorem~\ref{thm:zastava-comparison} identifies
$f_D(Z^{\c H}_M)$ with $Z_{(D,E)}$.
Hence the Kodaira embedding factors as
\[
\c I(D,E)\hookrightarrow Z_{(D,E)}
\hookrightarrow\b P(\b V_{(D,E)}).
\]
These embeddings are compatible with interval
inclusions and induce an inclusion
\[
j_p:
\Gr_{T,p}
=
\varinjlim_{\lambda\leq\mu}
\c I(\lambda[p],\mu[p])
\hookrightarrow
\c G^\kappa_p.
\]
The definition of the Kodaira embedding gives
$j_p^*\c O_{\c G^\kappa_p}(1) \simeq\b L^\kappa|_{\Gr_{T,p}}.$

\begin{prop}
	\label{prop:grassmannian-lattice-vertex-algebra}
	The restriction along $j_p$ induces
	canonical isomorphisms
	\[
	\Gamma\bigl(\c G^\kappa_p,\c O(1)\bigr)
	\xrightarrow{\sim}
	\Gamma\bigl(\Gr_{T,p},\b L^\kappa\bigr)
	\simeq L_\kappa^*.
	\]
	Equivalently,
	\[
	\varinjlim_{\lambda\leq\mu}
	\Gamma\bigl(
	Z_{(\lambda[p],\mu[p])},\c O(1)
	\bigr)^*
	\simeq L_\kappa,
	\]
	These isomorphisms are compatible with
	factorizations.
\end{prop}

\begin{proof}
	For $D\leq E$, define $r_{D,E}$ to be restriction
	along the Kodaira embedding:
	\[
	r_{D,E}:
	\Gamma\bigl(Z_{(D,E)},\c O(1)\bigr)
	\longrightarrow
	\Gamma\bigl(\c I(D,E),\b L^\kappa\bigr)
	=
	\b V_{(D,E)}^*.
	\]
	We first show that this map is an isomorphism. 	
	Put $M=E-D$. The following theorem is proved in \cite[\S2.4.3]{mirkovic2021loop}:
	$$r_{0,M}:
	\Gamma\bigl(Z^{\c H}_M,\c O(1)\bigr)
	\xrightarrow{\sim}
	\Gamma\bigl(\c H_M,\b L^\kappa\bigr).$$
	
	Choose a lift $\wt t_D$ of $t_D$ to the torus central extension as in the proof
	of Theorem~\ref{thm:interval-translation}.
	This lift gives isomorphisms
	$\c H_M\xrightarrow{\simeq}\c I(D,E)$ and
	$\c V_M\xrightarrow{\simeq}\b V_{(D,E)}$.
	Its projectivization restricts to $Z^{\c H}_M\simeq Z_{(D,E)}$ by
	Theorem~\ref{thm:zastava-comparison}. Using the isomorphism $\Phi_{\wt t_D}$ in Theorem~\ref{thm:interval-translation}, we get the following commutative diagram
	\[
	\begin{tikzcd}
		\Gamma\bigl(Z_{(D,E)},\c O(1)\bigr)
		\arrow[r,"r_{D,E}"]
		\arrow[d,"\Phi_{\wt t_D}^*"']
		&
		\Gamma\bigl(\c I(D,E),\b L^\kappa\bigr)
		\arrow[d,"\Phi_{\wt t_D}^*"]
		\\
		\Gamma\bigl(Z^{\c H}_M,\c O(1)\bigr)
		\arrow[r,"r_{0,M}"']
		&
		\Gamma\bigl(\c H_M,\b L^\kappa\bigr).
	\end{tikzcd}
	\]
	Hence $r_{D,E}$ is an isomorphism. 	One can further check the maps $r_{D,E}$ commute with inclusion of intervals. Taking inverse limits over intervals supported at $p$ we get
	\[	\Gamma\bigl(\c G^\kappa_p,\c O(1)\bigr)
		=\varprojlim_{\lambda\leq\mu}
		\Gamma\bigl(
		Z_{(\lambda[p],\mu[p])},\c O(1)
		\bigr)
	\simeq
	\varprojlim_{\lambda\leq\mu}
		\Gamma\bigl(
		\c I(\lambda[p],\mu[p]),\b L^\kappa
		\bigr)
		=\Gamma\bigl(\Gr_{T,p},\b L^\kappa\bigr)\simeq L_\kappa^*.
	\]
\end{proof}

\section{The Tits group action on $\c G^\kappa$}
\label{sec:tits-revised}
In this section we construct a Tits group action on the affine Grassmannian  $\c G^\kappa$. As before, fix an open embedding $C\hookrightarrow\b A^1_z$.
When $\kappa$ is of ADE type, the affine Grassmannian $\Gr_G$ doesn't carry a Weyl group action in general. A Tits subgroup of normalizer group $N_G(T)$ acts on $\Gr_G$. We want to extend the Tits group action for any symmetric simply-laced Kac--Moody Cartan matrix $\kappa$.

\subsection{The Tits group}
Suppose $\kappa$ is a symmetric simply-laced Kac--Moody Cartan matrix on $L$. Denote the simple coroots by $\{\alpha_i^\vee\}$, and identify $L=\b Z[I]$ with $\bigoplus_{i\in I}\b Z\alpha_i^\vee$. We have
$\kappa(\alpha_i^\vee,\alpha_i^\vee)=2$, while
$\kappa(\alpha_i^\vee,\alpha_j^\vee)\in\{0,-1\}$ for $i\ne j$.
The simple reflections $s_i(\lambda)=\lambda-
\kappa(\alpha_i^\vee,\lambda)\alpha_i^\vee$
generate the Coxeter--Weyl group $W_\kappa$.

We now define the Tits group directly from the Cartan datum. This definition does not use the corresponding Kac--Moody groups. Further discussions can be found in \cite[\S1.2]{ganapathy2024tits}.

\begin{defn}
	Put $T_2=L/2L$, and write $m_\lambda$ for the image of
	$\lambda\in L$ in $T_2$.  Thus
	$m_{\lambda+\mu}=m_\lambda m_\mu$ and $m_\lambda^2=1$.
	The \emph{Tits group} $\wt W_\kappa$ is generated by $T_2$ and elements
	$\wt s_i$, $i\in I$, satisfying
$$\wt s_i^2=m_{\alpha_i^\vee},\quad \wt s_i m_\lambda\wt s_i^{-1}=m_{s_i\lambda}$$
together with the braid relations
\[
\begin{aligned}
	\wt s_i\wt s_j=\wt s_j\wt s_i
	&&\text{if }\kappa(\alpha_i^\vee,\alpha_j^\vee)=0,
	\qquad
	\wt s_i\wt s_j\wt s_i=\wt s_j\wt s_i\wt s_j
	&&\text{if }\kappa(\alpha_i^\vee,\alpha_j^\vee)=-1.
\end{aligned}
\]
\end{defn}

There is an exact sequence (by definition it is not a central extension in general)
$	1\longrightarrow T_2\longrightarrow\wt W_\kappa
\longrightarrow W_\kappa\longrightarrow1,$
sending $\wt s_i$ to $s_i$.  
If $w=s_{i_1}\cdots s_{i_r}$ is reduced, set
$\wt w=\wt s_{i_1}\cdots\wt s_{i_r}.$
The braid relations imply that $\wt w$ is independent of the reduced
expression.  This defines a section (not a homomorphism) $W_\kappa\to\wt W_\kappa$.

\subsection{Intuition: group case}
\label{subsec:revised-group-case}

Assume temporarily that $\kappa$ is of finite ADE type, and let $G$ be
the corresponding simply-connected group.

For a  simple root $\alpha\in\Delta^+$, consider the corresponding $\sl_2$-triple $(e_\alpha\in \f g_\alpha, f_\alpha\in \f g_{-\alpha},h_\alpha\in \f{t})$. In  the corresponding  $\SL_2$-subgroup, we can take $\wt s_i$ as the order 4 element 
$\wt s_\alpha:=\exp(e_\alpha)\exp(-f_\alpha)\exp(e_\alpha)\in N_G(T)$.
Then $m_{\alpha_i^\vee}$ is realized as $\alpha_i^\vee(-1)$.

The Tits subgroup $\wt W_\kappa\subset G$ acts on $\Gr_G$.  For $\lambda\in X_*(T)$ write
$L_\lambda=z^\lambda\GO/\GO$.
For $w\in W_\kappa$, define
$$S^w_\lambda=
\left({w}N_\c K{w}^{-1}\right) z^\lambda\GO/\GO,\quad T^w_\lambda=
\left({w}N^-_\c K{w}^{-1}\right) z^\lambda\GO/\GO$$
Thus $S_\lambda^w$ and $T_\lambda^w$ depends only on $w$, because different lifts give the same groups. 

The lift is needed to define an actual map between semi-infinite orbits points in semiinfinite orbits. More precisely, if
$n\in N_{\c K}$, then
\[
\wt w\bigl(nL_\lambda\bigr)
=(\wt wn\wt w^{-1})L_{w\lambda}.
\]
It follows that
$	\wt wL_\lambda=L_{w\lambda},\ 
\wt wS_\lambda=S_{w\lambda}^w,\ 
\wt wT_\lambda=T_{w\lambda}^w.$

From the usual closure relations (see \cite[\S3]{mirkovic2007geometric})
\[
\ol S_{\lambda}
=\bigcup_{\mu\leq\lambda}S_{\mu},
\qquad
\ol T_{\lambda}
=\bigcup_{\mu\geq\lambda}T_{\mu}.
\] the action of Weyl group gives
\[
\ol S_{w\lambda}^w
=\bigcup_{\mu\leq\lambda}S_{w\mu}^w,
\qquad
\ol T_{w\lambda}^w
=\bigcup_{\mu\geq\lambda}T_{w\mu}^w.
\]
Therefore $w$ sends the boundary of $S_\lambda$ to the boundary of
$S_{w\lambda}^w$, and similarly for $T_\lambda$.

\subsubsection{Regular zastava fiber}

Let $D=\lambda[0]\leq E=\mu[0]$.  The comparison between interval
Zastava and the intersection of semi-infinite closures proved earlier
identifies $Z_{(\lambda[0],\mu[0])}$
with $\ol S_{-\lambda}\cap\ol T_{-\mu}$. Its action by $w$ can be defined so that $ wZ_{(\lambda[0],\mu[0])}
\simeq
\ol S_{-w\lambda}^w\cap
\ol T_{-w\mu}^w
.$

Suppose now that $\mu-\lambda=\alpha_j^\vee$ is a simple coroot. The zastava space $Z_{(\lambda[0],\mu[0])}$ is $\b P\bigl(\b L_\lambda^*\oplus\b L_\mu^*\bigr)$, which is isomorphic to $\b P^1$. It has two distinguished points  $\b P\b L_\lambda^*$ and
$\b P\b L_\mu^*$. 
By \cite[Theorem~3.1.2]{finkelberg2023elliptic}, they correspond to $L_{-\lambda},L_{-\mu}\in \Gr_G$ respectively.\footnote{The signs come from our Abel–Jacobi convention \(D\mapsto\mathcal O(-D)\).}

Since Weyl group $W$ preserves the bilinear form $\kappa$,   $\b L^\kappa$ has a $W$-equivariant structure. Hence 
$
\mathbb L_{\lambda[0]}^\vee
\longrightarrow
\mathbb L_{w\lambda[0]}^\vee, \mathbb L_{\mu[0]}^\vee
\longrightarrow
\mathbb L_{w\mu[0]}^\vee.
$
Consequently
$ w\,
\b P\bigl(\b L_\lambda^*\oplus\b L_\mu^*\bigr)
=\b P\bigl(\b L_{w\lambda}^*\oplus\b L_{w\mu}^*\bigr)$ is the image $wZ_{(\lambda[0],\mu[0])}$. It is the projective line joining the transformed marked
points. Note it need not to be the standard interval zastava 
$Z_{(w\lambda[0],w\mu[0])}$. If \(w\alpha_j^\vee\) is negative, the transformed endpoints do not define an interval. If it is positive but not simple, the transformed interval Zastava space can strictly contain this projective line.

\subsection{The Tits action on the Faltings line bundle}
\label{subsec:tits-faltings}
We have identified the lattice $L$ with the coroot lattice for $\kappa$, sending each basis $e_i$ to simple coroot $\alpha_i^\vee$.
Hence the Weyl group $W$ acts on $\Div_{C\times I}=\Div_C\otimes_\b Z L$.

\subsubsection{Tits group action on $\wh L$}
Let \(\widehat L\) be the central extension of $L$ by $\b G_m$ defined by the cocycle fixed in \S\ref{wh L}. Set $a_i=(1,e_i)\in\wh L$. 
The group $\wh L$ is generated by its central $\b G_m$ and
the $a_i$, with c onditions
$a_i a_j=(-1)^{\kappa(e_i,e_j)}a_j a_i$.

Let $\pi:\wh L\to L$ be the projection, and for any $x\in \wh L$, define $x_i:=\kappa(e_i,\pi(x))\in \b Z$.

\begin{lem}
	\label{lem:tits-lattice-lift}
	The maps $R_i: \wh L\to \wh L$, $R_i(x)=(-1)^{x_i(x_i-1)/2}x a_i^{-x_i}$ are group automorphisms, they define a $\wt W_\kappa$-action on $\wh L$ by sending $\wt s_i$ to $R_i$.
\end{lem}

\begin{proof}
	One can check that 
$$		R_i(a_j)=
		\begin{cases}
			-a_i^{-1},&j=i,\\
			a_j,&j\ne i,\ \kappa(e_i,e_j)=0,\\
			a_i a_j,&\kappa(e_i,e_j)=-1.
		\end{cases}$$
	Furthermore, for any $x$ we have $R_i^2(x)=(-1)^{x_i} x$.
	The properties follow by direct calculation, which we omit.
\end{proof}

\subsubsection{Tits group action on $\b L$}
We first recall the $\theta$-datum of the Faltings line bundle $\b L$ from Remark~\ref{thetadata}.
Let $j_\lambda:C\to\Div_{C\times I}$ be the map sending $p$ to $\lambda[p]$
. Consider line bundles over $C$:
\[\c E_\lambda
=(j_{-\lambda}^*\b L)^\vee
\otimes\omega_C^{\otimes\frac12\kappa(\lambda,\lambda)}.
\]
These line bundles are canonically trivialized. Let \(u_i\in\Gamma(C,\c E_{e_i}^\times)\) denote the corresponding trivializing sections.
The multiplication maps are
$\c E_\lambda\otimes\c E_\mu\xrightarrow{\sim}\c E_{\lambda+\mu}, 
xy=(-1)^{\kappa(\lambda,\mu)}yx.$
The commutator pairing is the same as the one for $\wh L$, hence there is an isomorphism of groups over $C$: 
$$\wh L\times C
\simeq \coprod_{\lambda\in L}\c E_\lambda^\times, \quad a_i=(1,e_i)\mapsto u_i.$$  
Transporting the $\wt W_\kappa$-action on \(\widehat L\) through this isomorphism gives a $\wt W_\kappa$-action on
$\coprod_{\lambda\in L}\c E_\lambda^\times$.
Each $\wt w\in\wt W_\kappa$ with image $w\in W_\kappa$
 induces isomorphisms of line bundles
$\c E_\lambda\xrightarrow{\sim}\c E_{w\lambda}.$

Since $w$ preserves $\kappa$, the cotangent terms in
$\c E_{\lambda}$ and $\c E_{w\lambda}$ coincide.
Removing the cotangent terms and dualizing
gives
\[
\phi_{\wt w}:\quad j_{-w\lambda}^*\b L\xrightarrow{\sim}j_{-\lambda}^*\b L.
\]
Hence we have defined a $\wt W$-equivariant structure on $\b L$. One can prove that it is compatible with factorization structure.

\subsection{Maps between interval vector spaces}
\label{subsec:tits-interval-maps}
The Weyl group action on $\Div_{C\times I}$ does not
preserve the subspace of effective divisors $\c H_{C\times I}$. Thus the assignment
$(D,E)\mapsto(wD,wE)$ does not define an action on $\Int$,
and we do not obtain an action on the vector bundle $\b V$
over that base. We instead construct compatible linear
inclusions between suitable fibers of $\b V$.

The following is obvious.
\begin{lem}
	\label{lem:tits-interval-bounds}
	For any interval $D\leq E$, there exists colored divisor $F\le G$ depending algebraically on $D,E$, such that there is a closed immersion
	\[
	w\c I(D,E)\hookrightarrow\c I(F,G).\qed
	\]
\end{lem}

For $\wt w\in \wt W_\kappa$ with image $w\in W$, the $\wt W_\kappa$-equivariant structure on $\b L$ gives
$$\Gamma(\c I(F,G),\b L)
\twoheadrightarrow\Gamma(w\c I(D,E),\b L)
\xrightarrow{\Phi_{\wt w}^*}\Gamma(\c I(D,E),\b L).$$
Dualizing yields
$A_{\wt w}^{D,E}:\b V_{(D,E)}\hookrightarrow\b V_{(F,G)}, 
a_{\wt w}^{D,E}=\b P(A_{\wt w}^{D,E}): \b P\b V_{(D,E)}\hookrightarrow\b P\b V_{(F,G)}.$
Assembling these maps gives a morphism $a_{\wt w}:\b P(\b V)\longrightarrow\b P(\b V)$.
This is a morphism of the total spaces, not an automorphism of the bundle over $\Int$.

Enlarging the target interval composes
$A_{\wt w}^{D,E}$ with the usual transition inclusion,
because restriction of sections is transitive.
Enlarging the source gives the corresponding commutative
diagram after choosing a common target.
Consequently different choices of bounds give the same
map after passage to a common larger interval.

The equivariant structure $\Phi$ also implies
$A_{\wt w}A_{\wt v}=A_{\wt w\wt v}, a_{\wt w}a_{\wt v}=a_{\wt w\wt v}$.
After placing the relevant maps in a common target interval, these maps give an action on the ambient
inductive limit of the projective interval spaces.

\subsection{Action on the affine Grassmannian $\c G^\kappa$}
\label{subsec:general-tits-action}
\subsubsection{Preservation of Zastava spaces}
We next prove that, after enlarging the target when needed,
the maps $a_{\wt w}$ carry Zastava spaces into Zastava spaces.
It suffices to treat the generators $\wt s_i$ of $\wt W_\kappa$.

We start with an interval $(0,e_j[p])$, one has
\[
Z_{(0,e_j[p])}
=\b P\bigl(\b L_0^\vee\oplus\b L_{e_j[p]}^\vee\bigr)
\simeq\b P^1.
\]
Its image under $a_{\wt s_i}$
is the projective line containing two points
labelled by $0$ and $s_i e_j[p]$ respectively.

\begin{lem}
	\label{lem:tits-elementary-line}
	The image of the projective line $Z_{(0,e_j[p])}$ satisfies
\[
a_{\wt s_i}^{0,e_j[p]}\bigl(Z_{(0,e_j[p])}\bigr)
\subset
\begin{cases}
	Z_{(-e_i[p],0)},&j=i,\\
	Z_{(0,e_j[p])},&j\ne i,\ \kappa(e_i,e_j)=0,\\
	Z_{(0,(e_i+e_j)[p])},&\kappa(e_i,e_j)=-1.
\end{cases}
\]
\end{lem}
\begin{proof}
	The three cases correspond to
	$s_i e_i=-e_i$, $s_i e_j=e_j$, and
	$s_i e_j=e_i+e_j$, respectively.
	Only the colors $i,j$ occur, hence by Lemma \ref{lem:based-product-decomposition} we can check with the
	corresponding group $G=SL_2$, $SL_2\times SL_2$,
	or $SL_3$. In the first two cases, the target is itself a projective
	line, and the image is the whole target.
	For the third case, the Tits group element acts the projective line to the line joining $L_0$ and $L_{-(e_i+e_j)}$ in $\Gr_{\SL_3}$.
	The curve lies in $\ol S_0\cap\ol T_{-(e_i+e_j)}
	\simeq Z_{(0,(e_i+e_j)[p])},$
	by Proposition~\ref{prop:interval-zastava-grassmannian}.
\end{proof}
We now apply the local projective space construction of
\S\ref{intervalzastava}. Over the left regular locus, the
Zastava fiber is a Segre embedding of product of
projective lines.
Since $a_{\wt s_i}$ is compatible with factorization,
Lemma~\ref{lem:tits-elementary-line} shows that it
sends this Segre product into the product of the
corresponding Zastava spaces, where these target need not be projective lines. By locality for Zastava spaces in \ref{cor:interval-locality}, their product is a Zastava fiber for
the sum of the target intervals. Hence for all left regular interval $(D,E)$, the image $a_{\wt s_i}(Z_{D,E})$ lies in some Zastava fiber $Z_{F,G}$. Taking closure we know the argument above holds for all interval $(D,E)$. Thus we obtain a $\wt W_\kappa$-action on
$\c G^\kappa$, compatible with factorization.

\begin{rmk}
	Let $\op{Out}(\kappa)$ be the automorphism group of the Dynkin diagram for $\kappa$. It acts on the affine Grassmannian $\c G^\kappa$ by	relabeling the colors. In ADE case, we know the automorphism group $O(\kappa,\b Z)$ preserving the bilinear form $\kappa$ is the semidirect product $W_\kappa\rtimes\op{Out}(\kappa)$ \cite[D.40]{fulton2013representation}.
\end{rmk}
\subsubsection{Comparison with the group case}

We compare the action constructed above with the usual
Tits action when $\kappa$ is of finite ADE type.
We first describe the action of torus elements
on $\b L$.

\begin{lem}
	\label{lem:torus-faltings-action}
	For $\nu,\lambda\in L\simeq X_*(T)$ and $c\in\b G_m$, the torus action on fibers of $\b L$ by
	\[
	\nu(c)\cdot v=c^{-\kappa(\nu,\lambda)}v,
	\qquad v\in\b L_{\lambda[p]}.
	\]
	As a corollary, the induced torus action on $\c E_\lambda$ acts by the same scalar.
\end{lem}

\begin{proof}
	Recall that the splitting $i:T_\c O\to\wh T_\kappa$
	gives a $T_\c O$-action on $\c R^\times=\wh T_\kappa/i(T_\c O)$
	by left multiplication. This action is trivial on the base	$\Gr_T=T_\c K/T_\c O$.
	Let $g\in T_\c K$, choose a lift
	$\widetilde g\in\wh T_\kappa$, and write
	$r=\widetilde g\,i(T_\c O)$ be the vector in $\c R_{gT_\c O}^\times$.
	For $u\in T_\c O$, we have
	$$	u\cdot r
	=i(u)\widetilde g\,i(T_\c O)
	=i(u)\widetilde g\,i(u)^{-1}\,i(T_\c O)
	=[i(u),\widetilde g]\,r.$$
	Thus the action on the fiber is given by the
	commutator pairing.	Taking $u=\nu(c)$ and $g=z^{-\lambda}$, the commutator
	formula in \S\ref{subsubsec:torus-central-extension}
	gives
	\[
	[i(\nu(c)),s(z^{-\lambda})]
	=\{c,z\}^{-\kappa(\nu,-\lambda)}
	=c^{\kappa(\nu,\lambda)}.
	\]
	Hence $\nu(c)$ acts on $\c R_{z^{-\lambda} T_\c O}$
	by this scalar.
	
	Using a local parameter at $p\in C$,
	Proposition~\ref{prop:faltings-heisenberg-comparison}
	gives $\b L_{\lambda[p]}
	\simeq\c R_{z^{-\lambda}T_\c O}^{\vee}.$
	Hence taking dual gives the asserted formula.
\end{proof}

We use $e_i,f_i$ for root vectors in the following
proposition, and write lattice elements using the
simple coroots $\alpha_i^\vee$.

\begin{prop}
	\label{prop:tits-group-comparison}
	Assume that $\kappa$ is of finite ADE type, and let
	$G$ be the corresponding simply-connected group.
	Keep the open embedding $C\hookrightarrow\b A^1_z$,
	and fix $p\in C$.
	There exist root vectors
	$e_i\in\mathfrak g_{\alpha_i}$ and
	$f_i\in\mathfrak g_{-\alpha_i}$, with
	$[e_i,f_i]=h_i$, such that, under the isomorphism
	$\c G^\kappa_p\simeq\Gr_{G,p}$, the action of
	$\wt s_i$ constructed in
	\S\ref{subsec:general-tits-action}
	agrees with the action of
	$n_i=\exp(e_i)\exp(-f_i)\exp(e_i)$
	for every $i\in I$.
\end{prop}

\begin{proof}
	By Proposition~\ref{prop:faltings-heisenberg-comparison}
	and \cite[Lemma~3.3.1(ii)]{zhu2009affine}, there are
	isomorphisms on $\Gr_{T,C^n}$
	\[
	\rho_n^*\b L
	\simeq \c R_{C^n}^{\vee}
	\simeq
	\left.\c O_{\Gr_{G,C^n}}(1)\right|_{\Gr_{T,C^n}},
	\]
	compatible with factorization.
	Here $\c O_{\Gr_{G,C^n}}(1)$ denotes the
	Beilinson--Drinfeld version of $\c O_{\Gr_G}(1)$.
	
	Start with any root vectors satisfying
	$[e_i,f_i]=h_i$, and put
	$n_i=\exp(e_i)\exp(-f_i)\exp(e_i)$.
	The usual action of $n_i$ on
	$\c O_{\Gr_{G,C^n}}(1)$ preserves its restriction
	to $\Gr_{T,C^n}$ and covers $s_i$.
	These actions are compatible with factorization.
	Through the preceding identifications, they induce
	isomorphisms
	\[
	\overline R_i:\c E_\lambda
	\xrightarrow{\sim}\c E_{s_i\lambda}
	\]
	preserving the multiplication maps.
	
	Recall $\wt s_i$ defines a map $R_i$ on
	$\coprod_{\lambda\in L}	\c E_\lambda^\times$, also sending 
	$\c E_\lambda$ to $\c E_{s_i\lambda}$. We now compare the maps $\overline R_i$ and $R_i$, their ratio on each $\c E_\lambda$ is an invertible regular function on $C$.
	However, the line bundles arise by restriction from the corresponding
	constructions on $\overline C=\b P^1$, see footnote \ref{footnote}.
	The ratio therefore belongs to
	$\Gamma(\overline C,\c O_{\overline C}^{\times})
	=\b C^\times$.
	Since both maps preserve multiplication, these
	constants define a character
	$\chi_i:L\to\b C^\times$, satisfying $	R_i(x)=\chi_i(\lambda)\overline R_i(x)$, for any $x\in\wh L$ whose image under \(\widehat L\to L\) is \(\lambda\).
	
	Suppose that $\kappa(\alpha_i^\vee,\lambda)=0$.
	The action of $z^\lambda$ on $e_i$ is $\Ad(z^\lambda)e_i=z^{\left\langle \alpha_i,\lambda \right\rangle }e_i=e_i$, hence $z^\lambda$ fixes the $\mathrm{SL}_2$-subgroup for triple $\{e_i,f_i,h_i\}$. Let $\varphi_i$ be embedding from the root $\SL_2$ to $G$, such that  $$\varphi_i\left( \begin{pmatrix}0&1\\-1&0\end{pmatrix}\right) =n_i.$$ 
	Because \(\mathcal O_{\mathrm{Gr}_G}(1)\) is \(G\)-equivariant, every \(h\in\mathrm{SL}_2\) gives a linear map $\mathcal O(1)_x
	\longrightarrow
	\mathcal O(1)_{\varphi_i(h)x}.$
	Now take \(x=z^\lambda G_{\mathcal O}\), then \(\varphi_i(h)x=x\) for every \(h\). So the maps above become linear automorphisms of the same fiber:
	\[
	\rho_x:\mathrm{SL}_2
	\to
	\operatorname{GL}\bigl(\mathcal O(1)_x\bigr)
	\simeq\mathbb G_m.
	\]
	Since \(\mathrm{SL}_2\) has no nontrivial characters, \(\rho_x\) is trivial. In particular, $n_i\cdot v=v$ for every $v\in\mathcal O(1)_x$.
		
	On the other hand, the fiber of \(\mathcal E_\lambda\) is
	\[
	(\mathcal E_\lambda)_p
	\simeq
	\bigl(\mathcal O_{\mathrm{Gr}_G}(1)_{z^\lambda G_{\mathcal O}}\bigr)^\vee
	\otimes
	(\omega_C)_p^{\otimes\frac12\kappa(\lambda,\lambda)}.
	\]
	The $\SL_2$-subgroup acts trivially on the fiber, thus \(\overline R_i\) acts as the identity. By the definition of $R_i$, their actions agree. Hence $\chi_i$ is trivial on the kernel of
	$\kappa(\alpha_i^\vee,-)$.
	There exists $c_i\in\b C^\times$ such that $\chi_i(\lambda)
	=c_i^{\kappa(\alpha_i^\vee,\lambda)}.$
	
	Replacing $e_i, f_i, n_i$ by $c_i e_i, c_i^{-1}f_i, \alpha_i^\vee(c_i)n_i$.
	By Lemma~\ref{lem:torus-faltings-action},
	the element $\alpha_i^\vee(c_i)$ acts on
	$\c E_{s_i\lambda}$ by
	$c_i^{-\kappa(\alpha_i^\vee,s_i\lambda)}
	=c_i^{\kappa(\alpha_i^\vee,\lambda)}.$
	The resulting group action therefore agrees with
	$R_i$ on the shifted $\theta$-datum. Hence $\{c_i e_i, c_i^{-1}f_i, h_i\}$ is the desired triple for all $i$.
\end{proof}

\bibliographystyle{alpha}
\bibliography{ref.bib}
\end{document}